\documentclass[12pt, runningheads]{amsart}
\usepackage{graphicx}
\font\bbr=msbm10
\font\sbbr=msbm7
\font\eightrm=cmr8
\font\gto=eufm10
\newcommand\ca[1]{{\mathcal #1}}
\newcommand\gt[1]{\hbox{\gto #1}}
\newcommand\longhookrightarrow{\lhook\joinrel\longrightarrow}
\newcommand\C{{\mathchoice{\hbox{\bbr C}}{\hbox{\bbr C}}{\hbox{\sbbr C}}
{\hbox{\sbbr C}}}}
\newcommand\N{{\mathchoice{\hbox{\bbr N}}{\hbox{\bbr N}}{\hbox{\sbbr N}}
{\hbox{\sbbr N}}}}
\newcommand\R{{\mathchoice{\hbox{\bbr R}}{\hbox{\bbr R}}{\hbox{\sbbr R}}
{\hbox{\sbbr R}}}}
\newcommand\Pro{{\mathchoice{\hbox{\bbr P}}{\hbox{\bbr P}}{\hbox{\sbbr P}}
{\hbox{\sbbr P}}}}
\newcommand\IM{\hbox{\rm Im}}
\newcommand\id{\hbox{\rm id}}
\newcommand\Hom{\hbox{\rm Hom}}
\newcommand\End{\hbox{\rm End}}
\newcommand\Sym{\hbox{\rm Sym}}
\newcommand\Sing{\hbox{\rm Sing}}
\let\de=\partial
\let\phe=\varphi
\newcommand{\void}{\ensuremath{\hbox{\bbr\symbol{"3F}}}}

\begin{document}

\title{Index Theorems for Holomorphic Maps and Foliations}
\thanks{Partially supported by Centro di Ricerca Matematica ``Ennio de Giorgi" and
the FIRB project \emph{Dinamica e azioni di gruppi su domini e
variet\`a complesse.}}

\author[M. Abate, F. Bracci, F. Tovena]{Marco Abate, Filippo Bracci \and Francesca Tovena }

\address{M. Abate: Dipartimento di Matematica, Universit\`a di Pisa, L.go Pontecorvo 5,
56127 Pisa, Italy} \email{abate@dm.unipi.it}

\address{F. Bracci \and F. Tovena: Dipartimento di Matematica, Universit\`a di Roma ``Tor Vergata'', Via
della Ricerca Scientifica 1, 00133 Roma, Italy}
\email{fbracci@mat.uniroma2.it, } \email{tovena@mat.uniroma2.it}

\emergencystretch15pt \frenchspacing

\newtheorem{theorem}{Theorem}[section]
\newtheorem{corollary}[theorem]{Corollary}
\newtheorem{proposition}[theorem]{Proposition}
\newtheorem{lemma}[theorem]{Lemma}

\theoremstyle{definition}\newtheorem{definition}[theorem]{Definition}

\theoremstyle{remark}
\newtheorem{remark}[theorem]{Remark}
\newtheorem{example}[theorem]{Example}
\newtheorem{case}{Case}

\begin{abstract}
We describe a general construction providing index theorems
localizing the Chern classes of the normal bundle of a subvariety
inside a complex manifold. As particular instances of our
construction we recover both Lehmann-Suwa's generalization of the
classical Camacho-Sad index theorem for holomorphic foliations and
our index theorem for holomorphic maps with positive dimensional
fixed point set. Furthermore, we also obtain generalizations of
recent index theorems of Camacho-Movasati-Sad and Camacho-Lehmann
for holomorphic foliations transversal to a subvariety.
\end{abstract}

\maketitle

\section{Introduction}
\label{intro}
In 1982, C. Camacho and P. Sad \cite{CS} proved the existence of separatrices for
singular holomorphic foliations in dimension~2. One of the main tools in
their proof was the following index theorem:

\begin{theorem}[Camacho-Sad]\label{th:CS}
Let $S$ be a compact Riemann surface embedded in a smooth complex surface~$M$. Let
$\ca F$ be a one-dimensional singular holomorphic foliation defined in a
neighbourhood of~$S$ and such that $S$ is a leaf of~$\ca F$, that is such that $\ca F$ is
tangent to~$S$. Then it is possible to associate to any singular point~$q\in S$
of~$\ca F$ a complex number~$\iota_q(\ca F,S)\in\C$,
the \emph{index} of the foliation along~$S$ at~$q$, depending only on the local
behavior of~$\ca F$ near~$q$, such that
\[
\sum_{q\in\hbox{\eightrm Sing}(\ca F)}\iota_q(\ca F,S)=\int_S c_1(N_S),
\]
where $c_1(N_S)$ is the first Chern class of the normal bundle~$N_S$ of~$S$ in~$M$.
\end{theorem}

The index $\iota_q(\ca F,S)$ can be explicitely computed using a local vector field
generating $\ca F$ in a neighbourhood of~$q$ (see Example~\ref{ex:Filippo2}).

Thus Theorem~\ref{th:CS} gives a quantitative connection between the way $S$ sits in~$M$
(the integral of the first Chern class of~$N_S$ is equal to the self-intersection
number~$S\cdot S$ of~$S$) and the behavior of singular foliations tangent to~$S$.

Due to its importance (see, e.g., \cite{Li} and \cite{Bru} for applications), in the
next twenty years this theorem has been generalized in several ways; see, e.g., Lins Neto
\cite{Li}, Suwa \cite{S1}, Lehmann~\cite{Le}, Lehmann-Suwa~\cite{LS1}, \cite{LS2},
and references therein (see also~\cite{BS1}, \cite{BS2}, where also the ambient~ space~$M$ is
allowed to be singular). In particular, using
\v Cech-de Rham cohomology, Lehmann and Suwa (see
\cite{S2} for a systematic exposition) proved what we are going to consider as a model index
theorem: if the possibly singular holomorphic foliation~$\ca F$ of dimension~$\ell$ is tangent to a
possibly singular (but not too wild: see Definition~\ref{def:nest} and Example~\ref{ex:nest})
subvariety~$S$ of dimension~$d$ in the ambient manifold~$M$, then the Chern classes of the normal
bundle~$N_S$ of degree higher than $d-\ell$ can be localized at the singularities --- that is,
obtained as sum of local residues depending only on the behavior of~$S$ and~$\ca F$ nearby the
singularities (of $S$ and of $\ca F$ in~$S$). We explicitely remark that one of the main
ingredients in their proof is the construction of a partial holomorphic connection on the
the normal bundle~$N_S$ outside a suitable analytic subset of~$S$ (where ``partial" here means that
we are differentiating only along some tangent directions, the ones contained in~$\ca F$).

The results obtained in those papers are apparently strictly inside the theory of holomorphic
foliations; the arguments used needed the existence of the foliation~$\ca F$ in a
neighbourhood of the subvariety~$S$, and the tangency of~$\ca F$ to~$S$. In the last five
years, however, a number of results have appeared strongly suggesting that these might be
unnecessary limitations: the tangency of $\ca F$ to~$S$ might be replaced by hypotheses on the
embedding of~$S$ into the ambient space~$M$ (\cite{C}, \cite{CL}, \cite{CMS}), and, perhaps more
strikingly, the foliation can be replaced by a holomorphic self-map of the ambient manifold fixing
pointwise the subvariety~$S$ (\cite{A}, \cite{BT}, \cite{ABT}). Furthermore, there was the
tantalizing fact that the statements of all these new index theorems were clearly
similar, and yet they all needed slightly different proofs; none of them was consequence of any of
the others.

The main goal of this paper is to show how it is possible to recover all these
index theorems (and a couple of new ones) using a universal construction
having a priori nothing to do with either foliations or self-maps. More precisely, we
shall reduce the proof of such an index theorem to the construction of an $\ca
O_{S^o}$-morphism $\psi\colon\ca F\to\ca A$ satisfying a splitting condition (see
Theorem~\ref{th:sommario}.(ii) for the exact condition), where:
$S^o$ is the complement in~$S$ of the singular points of~$S$ and of the singular points of the
object (foliation or self-map) we are interested in; $\ca F$ is the sheaf of germs of holomorphic
sections of a suitable sub-bundle of the tangent bundle~$TS^o$; and $\ca A$ is an universal $\ca
O_{S^o}$-locally free sheaf depending only on the embedding of~$S^o$ into the ambient space~$M$.
The details of the construction of~$\psi$ will of course depend on the particular situation we are
dealing with (foliation or self-map, tangential or transversal); but as soon as such a $\psi$
exists then an index theorem analogous to Lehmann-Suwa's model one follows (see
Theorem~\ref{th:generalindex}). We shall also show (in Sections~7 and~8) how to construct such
a~$\psi$ in several cases, and we shall be able to recover all the index theorems of this kind
known up to now (for $M$ smooth), together with two new ones: the first one
(Theorem~\ref{th:foltrans}) on foliations transverse to~$S$ generalizes both Camacho-Movasati-Sad's
(\cite{CMS}) and Camacho-Lehmann's (\cite{C},
\cite{CL}) results, while the second one (Theorem~\ref{th:maphcd}) extends
the results of~\cite{ABT}. We also explicitely compute the index at isolated singularities in a
simple but important case using a Grothendieck residue; see Remark~\ref{rem:Filippo1}
and Example~\ref{ex:Filippo2}.

Very briefly, the reason why the existence of a morphism $\psi\colon\ca F\to\ca A$ implies an
index theorem can be explained as follows. The sheaf $\ca A$ comes provided with an important
additional structure: an universal holomorphic
$\theta_1$-connection (see Sections~4 and~5 for precise definitions) on the normal bundle~$N_{S^o}$
of~$S^o$ in~$M$. Then, following ideas due to Atiyah~\cite{At}, we shall be able to prove that the
existence of such a morphism~$\psi$ is equivalent to the existence of a partial
holomorphic connection on~$\ca N_{S^o}$ along~$\ca F$. Having this, an
argument essentially due to Baum and Bott (\cite{BB}; see Theorem~\ref{th:Bott}) yields the
vanishing on~$S^o$ of suitable Chern classes of~$N_{S^o}$; and then a general cohomological argument
(developed by Lehmann and Suwa \cite{Le},
\cite{S2}) allows one to infer from this vanishing the localization at the
singularities of the corresponding Chern classes --- that is, an index
theorem.

Let us finally describe the plan of the paper. In Sections~2 and~3 we collect a number of
definitions and properties concerning infinitesimal neighbourhoods of subvarieties that we shall
need in the rest of the paper, and we describe the conditions we shall impose on the embedding of
the subvariety into the ambient manifold to deal with the transversal cases (we
also refer to~\cite{ABT3} for more details on these conditions). In Sections~4 and~5 we introduce
the sheaf~$\ca A$ and its additional structures, while in Section~6 we prove the vanishing
Theorem~\ref{th:Bott} and our general index Theorem~\ref{th:generalindex}.
In Section~7 we show how to build the morphism~$\psi$ for holomorphic foliations (and, in
particular, we get the new Theorem~\ref{th:foltrans}), and finally in Section~8 we do the same for
holomorphic self-maps (and, in particular, we get the new Theorem~\ref{th:maphcd}).

\bigskip

We would like to thank Francesco Russo for
pointing out reference \cite{MP}, Jorge V. Pereira and Tatsuo Suwa
for some useful conversations, and the Departments of Mathematics of Jagellonian
University, Krak\'ow, and of University of Michigan, Ann Arbor, for the warm hospitality offered
to the first and third author during the preparation of this paper.

\section{Splitting submanifolds}
\label{sec:2}

Let us begin by recalling a few general facts on sequences of sheaves.
We say that an exact sequence of sheaves (of abelian groups, rings,
modules\dots)
\[
O\longrightarrow\ca R\mathop{\longrightarrow}\limits^\iota\ca S
\mathop{\longrightarrow}\limits^p \ca T\longrightarrow O
\]
on a variety $S$ \emph{splits} if there is a morphism $\sigma\colon\ca T\to\ca S$ of sheaves
(of abelian groups, rings, modules\dots) such that $p\circ\sigma=\id$. Any such morphism is
called a \emph{splitting morphism.} A morphism of sheaves of abelian groups $\tau\colon\ca
S\to\ca R$ such that $\tau\circ\iota=\id$ is called a \emph{left splitting morphism.}

The following facts are well-known, and easy to prove:

\begin{lemma}
\label{th:uquattro}
Let
\begin{equation}
O\longrightarrow\ca R\mathop{\longrightarrow}\limits^\iota\ca S
\mathop{\longrightarrow}\limits^p \ca T\longrightarrow O
\label{eq:seqr}
\end{equation}
be an exact sequence of sheaves of abelian groups over a variety~$S$. Then:
\begin{itemize}
\item[\textup{(i)}]the sequence $(\ref{eq:seqr})$
splits if and only if there exists a left splitting morphism $\tau\colon\ca S\to\ca R$;
\item[\textup{(ii)}]if $(\ref{eq:seqr})$ splits, for any splitting morphism
$\sigma\colon\ca T\to\ca S$ there exists a unique left splitting morphism $\tau\colon\ca
S\to\ca R$ such that
$\tau\circ\sigma=O$ and
\[
\iota\circ\tau+\sigma\circ p=\id;
\]
\item[\textup{(iii)}]if $(\ref{eq:seqr})$ splits, then there is
a $1$-to-$1$ correspondance between splitting morphisms and elements of
$H^0\bigl(S,\Hom(\ca
T,\ca R)\bigr)$. More precisely, if $\sigma_0\colon\ca T\to\ca
S$ is a splitting morphism, all others splitting morphisms are of the
form~$\sigma_0-\iota\circ\phe$ with $\phe\in H^0\bigl(S,\Hom(\ca T,\ca R)\bigr)$, while if
$\tau_0\colon\ca T\to\ca S$ is a left splitting morphism, all others left splitting morphisms
are of the form
$\tau_0+\phe\circ p$ with~$\phe\in H^0\bigl(S,\Hom(\ca
T,\ca R)\bigr)$.
\end{itemize}
\end{lemma}

Following Grothendieck and Atiyah, we can give an useful cohomological characterization of
splitting for sequences of locally free $\ca O_S$-modules. Let
\begin{equation}
O\longrightarrow\ca E'\longrightarrow\ca E\longrightarrow\ca E''\longrightarrow O
\label{eq:seqE}
\end{equation}
be an exact sequence of sheaves of locally free $\ca O_S$-modules. Applying the functor
$\Hom(\ca E'',\cdot)$ to this sequence we get the exact sequence
\begin{equation}
O\longrightarrow\Hom(\ca E'',\ca E')\longrightarrow\Hom(\ca E'',\ca E)\longrightarrow
\Hom(\ca E'',\ca E'')\longrightarrow O.
\label{eq:seqEH}
\end{equation}
Let $\delta\colon H^0\bigl(S,\Hom(\ca E'',\ca E'')\bigr)\to H^1\bigl(S,\Hom(\ca E'',\ca
E')\bigr)$ be the connecting homomorphism in the long exact cohomology sequence
of~(\ref{eq:seqEH}). Then we can associate to the exact sequence (\ref{eq:seqE}) the cohomology
class
\[
\delta(\id_{\ca E''})\in H^1\bigl(S,\Hom(\ca E'',\ca E')\bigr).
\]

This procedure gives a $1$-to-$1$ correspondance between the
  group  $H^1\bigl(S,\Hom(\ca E'',\ca E')\bigr)$ and
isomorphism classes of exact sequences of locally free $\ca
O_S$-modules starting with~$\ca E'$ and ending with~$\ca E''$.
Indeed, we have (see \cite{At}, Proposition~1.2):

\begin{proposition}
\label{th:uotto}
Let $S$ be a complex manifold. Then two exact sequences of locally free $\ca O_S$-modules are
isomorphic if and only if they correspond to the same cohomology class. In particular, an exact
sequence~$(\ref{eq:seqE})$ of locally free $\ca O_S$-modules splits if and only if it corresponds to the zero
cohomology class.
\end{proposition}

Let us now introduce the sheaves (and sequences of sheaves) we are interested in. Let $M$ be
a complex manifold of dimension $n$, and let $S$ be a reduced, globally irreducible subvariety
of~$M$ of codimension~$m\ge 1$. We denote: by $\ca O_M$ the sheaf of germs of holomorphic
functions on~$M$; by~$\ca I_S$ the subsheaf of~$\ca O_M$ of germs vanishing on~$S$; and by~$\ca
O_S$ the quotient sheaf~$\ca O_M/\ca I_S$ of germs of holomorphic functions on~$S$.
Furthermore, let~$\ca T_M$ denote the sheaf of germs of holomorphic sections of the
holomorphic tangent bundle~$TM$ of~$M$, and $\Omega_M$ the sheaf of germs
of holomorphic 1-forms on~$M$. Finally, we shall denote by~$\ca T_{M,S}$ the sheaf of
germs of holomorphic sections along~$S$ of the restriction~$TM|_S$ of~$TM$ to~$S$, and by
$\Omega_{M,S}$ the sheaf of germs of holomorphic sections along~$S$ of~$T^*M|_S$. It is
easy to check that $\ca T_{M,S}=\ca T_M\otimes_{\ca O_M}\ca O_S$
and~$\Omega_{M,S}=\Omega_M\otimes_{\ca O_M}\ca O_S$.

For $k\ge 1$ we shall denote by $f\mapsto[f]_k$ the canonical projection of~$\ca O_M$ onto~$\ca
O_M/\ca I_S^k$. The
\emph{cotangent sheaf}~$\Omega_S$ of~$S$ is defined by
\[
\Omega_S=\Omega_{M,S}\big/d_2(\ca I_S/\ca I_S^2),
\]
where $d_2\colon\ca O_M/\ca I_S^2\to\Omega_{M,S}$ is given by $d_2[f]_2=df\otimes[1]_1$.
In particular, we have the \emph{conormal sequence} of sheaves of~$\ca O_S$-modules associated
to~$S$:
\[
\ca I_S/\ca I_S^2\mathop{\longrightarrow}\limits^{d_2}\Omega_{M,S}
\mathop{\longrightarrow}\limits^{p}\Omega_S\longrightarrow O.
\]
Applying the functor $\Hom_{\ca O_S}(\cdot,\ca O_S)$ to the conormal sequence we get the
\emph{normal sequence} of sheaves of $\ca O_S$-modules associated to~$S$:
\[
O\longrightarrow \ca T_S \longhookrightarrow\ca T_{M,S}
\mathop{\longrightarrow}\limits^{p_2}\ca N_S,
\]
where $\ca T_S=\Hom_{\ca O_S}(\Omega_S,\ca O_S)$ is the \emph{tangent sheaf} of~$S$, $p_2$
is the morphism dual to~$d_2$, and $\ca
N_S=\Hom_{\ca O_S}(\ca I_S/\ca I_S^2,\ca O_S)$ is the \emph{normal sheaf} of~$S$.

As mentioned in the introduction, to get index theorems in the transversal case we shall need
hypotheses on the embedding of~$S$ into~$M$. The first such hypothesis is:

\begin{definition}
Let $S$ be a reduced, globally irreducible subvariety of a complex
manifold $M$. We say
that $S$ {\sl splits} into~$M$ if there exists a morphism of sheaves of
$\ca O_S$-modules
$\sigma\colon\Omega_S\to\Omega_{M,S}$ such that $p\circ\sigma=\id$, where
$p\colon\Omega_{M,S}\to\Omega_S$ is the canonical projection.
\end{definition}

\begin{remark}
\label{rem:smooth}
It is not difficult to prove that if $S$ splits into $M$ then it is
necessarily non-singular, and the morphism $d_2\colon\ca I_S/\ca I_S^2\to\Omega_{M,S}$ is
injective. In particular, when $S$ splits into~$M$ the sequence
\begin{equation}
O\longrightarrow\ca I_S/\ca I_S^2\mathop{\longrightarrow}\limits^{d_2}\Omega_{M,S}
\mathop{\longrightarrow}\limits^{p}\Omega_S\longrightarrow O
\label{eq:sequnoprimo}
\end{equation}
is a splitting exact sequence of locally free $\ca O_S$-modules, and we also have a left splitting
morphism $\tau\colon\Omega_{M,S}\to\ca I_S/\ca I_S^2$.
\end{remark}

We shall now describe several equivalent characterizations
of splitting subvarieties. In doing so, we shall introduce notations and terminologies that
shall be useful in the rest of the paper.

\begin{definition}
Let $S$ be a reduced, globally irreducible subvariety of a
complex manifold $M$. For any~$k\ge 1$ let $\theta_k\colon\ca O_M/\ca I_S^{k+1}\to\ca
O_M/\ca I_S$ be the canonical projection given by $\theta_k([f]_{k+1})=[f]_1$. The
\emph{$k$-th infinitesimal neighbourhood} of $S$ in $M$ is the ringed space $S(k)=(S,\ca
O_M/\ca I_S^{k+1})$ together with the canonical inclusion of
ringed spaces $\iota_k\colon S=S(0)\to S(k)$ given by~$\iota_k=(\id_S,\theta_k)$. We also
put~$\ca O_{S(k)}=\ca O_M/\ca I_S^{k+1}$. A \emph{$k$-th
order lifting} is a splitting morphism
$\rho\colon\ca O_S\to\ca O_{S(k)}$ for the exact sequence of sheaves of rings
\[
O\longrightarrow\ca I_S/\ca I_S^{k+1}\longhookrightarrow\ca
O_{S(k)}\mathop{\longrightarrow}\limits^{\theta_k}\ca O_S\longrightarrow O.
\]
\end{definition}

\begin{definition}
Let $\ca O$, $\ca R$ be sheaves of rings, $\theta\colon\ca R\to\ca O$ a
morphism of sheaves of rings, and $\ca M$ a sheaf of $\ca O$-modules. A {\sl
$\theta$-derivation} of~$\ca R$ in~$\ca M$ is a morphism of sheaves of abelian groups
$D\colon\ca R\to\ca M$ such that
$$
D(r_1r_2)=\theta(r_1)\cdot D(r_2)+\theta(r_2)\cdot D(r_1)
$$
for any $r_1$, $r_2\in\ca R$. In other words, $D$ is a
derivation with respect to the $\ca R$-module structure induced via
restriction of scalars by~$\theta$.
\end{definition}

We can now give a first list of conditions equivalent to splitting (see \cite{G}, p. 373,
\cite{MP}, Lemma~1.1, and \cite{Ei}, Proposition~16.12 for proofs):

\begin{proposition}
\label{th:ucinque}
Let $S$ be a reduced, globally irreducible
subvariety of a complex manifold~$M$. Then there is a $1$-to-$1$
correspondance among the following classes of morphisms:
\begin{itemize}
\item[\textup{(a)}] morphisms $\sigma\colon\Omega_S\to\Omega_{M,S}$ of sheaves of
$\ca O_S$-modules such that $p\circ\sigma=\id$;
\item[\textup{(b)}] morphisms $\tau\colon\Omega_{M,S}\to\ca I_S/\ca I_S^2$ of sheaves of $\ca
O_S$-modules such that $\tau\circ d_2=\id$;
\item[\textup{(c)}] $\theta_1$-derivations $\tilde\rho\colon\ca O_{S(1)}\to\ca I_S/\ca I_S^2$
such that $\tilde\rho\circ i_1=\id$, where
$i_1\colon\ca I_S/\ca I_S^2\hookrightarrow\ca O_{S(1)}$ is the
canonical inclusion;
\item[\textup{(d)}] morphisms $\rho\colon\ca O_S\to\ca O_M/\ca I_S^2$ of sheaves of rings such
that
$\theta_1\circ\rho=\id$.
\end{itemize}
In particular,
$S$ splits into $M$ if and only if it has a first order lifting. Finally, if any (and hence all) of
the classes \textup{(a)--(d)} is not empty, then it is in $1$-to-$1$ correspondance with the
following classes of morphisms:
\begin{itemize}
\item[\textup{(e)}] morphisms $\tau^*\colon\ca N_S\to\ca T_{M,S}$ of sheaves of $\ca O_S$-modules
such that $p_2\circ\tau^*=\id$;
\item[\textup{(f)}] morphisms $\sigma^*\colon\ca T_{M,S}\to\ca T_S$ of sheaves of
$\ca O_S$-modules such that $\iota\circ\sigma^*=\id$, where $\iota\colon\ca T_S\to\ca T_{M,S}$ is
the canonical inclusion.
\end{itemize}
\end{proposition}

We have already noticed (Remark~\ref{rem:smooth}) that a splitting subvariety is necessarily
non-singular; therefore we can use differential geometric techniques to get another couple of
characterizations of splitting submanifolds.

\begin{definition}
Let $S$ be a (not necessarily closed) complex submanifold of codimension $m\ge1$ in a
complex manifold $M$ of dimension~$n\ge 2$, and let~$(U_\alpha,z_\alpha)$ be a chart of $M$. We
shall sistematically write
$z_\alpha=(z_\alpha^1,\ldots,z_\alpha^n)=(z'_\alpha,z''_\alpha)$, with
$z'_\alpha=(z_\alpha^1,\ldots,z_\alpha^m)$ and
$z''_\alpha=(z_\alpha^{m+1},\ldots,z_\alpha^n)$. We shall say that $(U_\alpha,z_\alpha)$ is
\emph{adapted} to~$S$ if either $U_\alpha\cap S=\void$ or
\[
U_\alpha\cap S=\{z_\alpha^1=\cdots=z_\alpha^m=0\}.
\]
In particular, if
$(U_\alpha,z_\alpha)$ is adapted to~$S$ then $\{z_\alpha^1,\ldots,z_\alpha^m\}$ is a set of
generators of~$\ca I_{S,x}$ for all~$x\in U_\alpha\cap S$. An atlas $\gt
U=\{(U_\alpha,z_\alpha)\}$ of~$M$ is \emph{adapted} to~$S$ if all its charts are; then $\gt
U_S=\{(U_\alpha\cap S,z''_\alpha)\mid U_\alpha\cap S\ne\void\}$ is an atlas for~$S$. The
\emph{normal bundle}~$N_S$ of~$S$ in~$M$ is the quotient bundle~$TM|_S/TS$; its dual is the
\emph{conormal bundle}~$N^*_S$. If
$(U_\alpha,z_\alpha)$ is a chart adapted to~$S$, for~$r=1,\ldots,m$ we shall denote
by~$\de_{r,\alpha}$ the projection of~$\de/\de z_\alpha^r|_{U_\alpha\cap S}$ in~$N_S$, and
by~$\omega_\alpha^r$ the local section of~$N_S^*$ induced by~$dz_\alpha^r|_{U_\alpha\cap
S}$. Then $\{\de_{1,\alpha},\ldots,\de_{m,\alpha}\}$
and~$\{\omega_\alpha^1,\ldots,\omega_\alpha^m\}$ are local frames over~$U_\alpha\cap S$ for
$N_S$ and~$N_S^*$ respectively, dual to each other.
\end{definition}

\begin{remark}
From now on, every chart and atlas we consider on~$M$ will be adapted to~$S$.
We shall use Einstein convention on the sum over repeated indices. Indices like
$j$, $h$, $k$ will run from~1 to~$n$; indices like $r$, $s$, $t$,~$u$, $v$ will run from~1
to~$m$; and indices like $p$, $q$ will run from~$m+1$ to~$n$.
\end{remark}

\begin{remark}
If $(U_\alpha,z_\alpha)$ and $(U_\beta,z_\beta)$ are two adapted charts with
$U_\alpha\cap U_\beta\cap S\ne\void$, then it is easy to check that
\[
\left.\frac{\de z_\beta^r}{\de z_\alpha^p}\right|_S\equiv O
\]
for all $r=1,\ldots,m$ and $p=m+1,\ldots,n$.
\end{remark}

Then computing the cohomology class associated to the conormal sequence~(\ref{eq:sequnoprimo}) and
recalling Proposition~\ref{th:uotto} we get:

\begin{proposition}
\label{th:unove}
Let $S$ be a complex submanifold of codimension $m$ of a
complex manifold~$M$, and let $\gt U=\{(U_\alpha,z_\alpha)\}$ be an adapted atlas. Then the
cohomology class $\gt s\in H^1\bigl(S,\Hom(\Omega_S,\ca N^*_S)\bigr)$ associated to the conormal
exact sequence of~$S$ is represented by the $1$-cocycle $\{\gt s_{\beta\alpha}\}\in H^1\bigl(\gt
U_S,\Hom(\Omega_S,\ca N^*_S)\bigr)$ given by
\[
\gt s_{\beta\alpha}=-\left.\frac{\de z_\beta^r}{\de z_\alpha^s}\,\frac{\de z_\alpha^p}{\de
z_\beta^r}\right|_S\omega_\alpha^s\otimes\frac{\de}{\de z_\alpha^p}\in H^0(U_\alpha\cap
U_\beta\cap S,\ca N_S^*\otimes\ca T_S).
\]
In particular, $S$ splits into $M$ if and only if $\gt s=O$.
\end{proposition}

We can rewrite this characterization in a more useful form using the notion of splitting
atlas, originally introduced in~\cite{ABT}.

\begin{definition}
Let $\gt U=\{(U_\alpha,z_\alpha)\}$ be an adapted atlas for a complex
submanifold~$S$ of codimension~$m\ge 1$ of a complex $n$-dimensional manifold~$M$. We say that
$\gt U$ is a \emph{splitting atlas} if
\[
\left.\frac{\de z_\beta^p}{\de z_\alpha^r}\right|_S\equiv O
\]
for all $r=1,\ldots, m$, $p=m+1,\ldots,n$ and indices $\alpha$, $\beta$ so that
$U_\alpha\cap U_\beta\cap S\ne\void$.
\end{definition}

\begin{definition}
Let $\gt U=\{(U_\alpha,z_\alpha)\}$ be an atlas adapted to~$S$. If $\rho\colon\ca O_S\to\ca
O_{S(1)}$ is a first order lifting for~$S$, we say~$\gt U$ is
\emph{adapted} to~$\rho$ if
\begin{equation}
\rho([f]_1)=[f]_2-\left[\frac{\de f}{\de z^r_\alpha} z_\alpha^r\right]_2
\label{eq:defgrho}
\end{equation}
for all $f\in\ca O(U_\alpha)$ and all indices $\alpha$ such that $U_\alpha\cap S\ne\void$.

\begin{remark}
In \cite{ABT3} it is shown that $\gt U$ is adapted to~$\rho$ if and only if
\[
\rho(g)(z_\alpha)=g_\alpha(O',z''_\alpha)+\ca I_S^2
\]
for all $g\in\ca O(U_\alpha\cap S)$ and all indices $\alpha$
such that~$U_\alpha\cap S\ne\void$, where we are assuming (without loss of generality) that
$z_\alpha\in U_\alpha$ implies $(O',z''_\alpha)\in U_\alpha\cap S$.
\end{remark}
\end{definition}

\begin{proposition}
\label{th:udieci}
Let $S$ be a complex submanifold of codimension~$m\ge 1$ of a
$n$-dimensional complex manifold~$M$. Then:
\begin{itemize}
\item[\textup{(i)}]$S$ splits into $M$ if and only if there exists a splitting atlas
for $S$ in $M$;
\item[\textup{(ii)}]an atlas adapted to $S$ is splitting if and only if
it is adapted to a first order lifting;
\item[\textup{(iii)}]if $S$ splits into $M$, then for any first order lifting there exists an
atlas adapted to it.
\end{itemize}
\end{proposition}

\begin{proof}
(i) By Propositions~\ref{th:uotto} and~\ref{th:unove}, the existence of a splitting atlas clearly
implies that $S$ splits into $M$. Conversely, assume that $S$ splits into $M$. Then by
Propositions~\ref{th:uotto} and~\ref{th:unove} we can find an adapted atlas $\gt U$ and a
0-cochain
$\gt c=\{\gt c_\alpha\}\in H^0(\gt U_S,\ca N_S^*\otimes\ca T_S)$ such that $\gt
s_{\beta\alpha}=\gt c_\beta-\gt c_\alpha$ on $U_\alpha\cap\ U_\beta\cap S$, that is
\begin{equation}
-\left.\frac{\de z_\alpha^r}{\de z_\beta^s}\,\frac{\de z_\beta^q}{\de z_\alpha^r}\right|_S
=(c_\alpha)^p_r \left.\frac{\de z_\alpha^r}{\de z_\beta^s}\,\frac{\de z_\beta^q}{\de
z_\alpha^p}\right|_S - (c_\beta)^q_s
\label{eq:coh}
\end{equation}
on $U_\alpha\cap U_\beta\cap S$ for all $s=1,\ldots,m$ and $p=m+1,\ldots,n$, where we have
written
\[
\gt c_\alpha=(c_\alpha)^p_s\,\omega_\alpha^s\otimes\frac{\de}{\de z_\alpha^p}.
\]
Then using (\ref{eq:coh}) it is easy to check that the coordinates
\begin{equation}\label{eq:ccuno}
\begin{cases}\hat z_\alpha^s=z_\alpha^s,\\
\hat{z}^p_\alpha=z_\alpha^p-(c_\alpha)^p_r z_\alpha^r, \end{cases}
\end{equation}
restricted to suitable open subsets $\hat U_\alpha\subseteq U_\alpha$, give a splitting
atlas~$\hat{\gt U}=\{(\hat U_\alpha,\hat z_\alpha)\}$.
\smallskip

(ii) Let $\gt U=\{(U_\alpha,z_\alpha)\}$ be an atlas adapted to~$S$. Setting
\[
\rho_\alpha([f]_1)=[f]_2-\left[\frac{\de f}{\de z^r_\alpha} z_\alpha^r\right]_2
\]
we define local first order
liftings $\rho_\alpha\colon\ca O_S|_{U_\alpha\cap S}\to\ca O_{S(1)}|_{U_\alpha\cap S}$; we
claim that in this way we get a global
first order lifting~$\rho$ if and only if $\gt U$ is a splitting atlas (necessarily adapted
to~$\rho$). But indeed
\begin{eqnarray*}
\rho_\beta([f]_1)-\rho_\alpha([f]_1)&=&\left[\frac{\de f}{\de z^r_\alpha}
z_\alpha^r\right]_2-\left[\frac{\de f}{\de z^s_\beta} z_\beta^s\right]_2\\
&=&\left[\frac{\de f}{\de z_\beta^s}\left(\frac{\de z^s_\beta}{\de
z^r_\alpha}z_\alpha^r-z^s_\beta\right)\right]_2+
\left[\frac{\de f}{\de z^p_\beta}\frac{\de z^p_\beta}{\de z^r_\alpha}z_\alpha^r
\right]_2\\
&=&\left[\frac{\de f}{\de z^p_\beta}\frac{\de z^p_\beta}{\de z^r_\alpha}z_\alpha^r
\right]_2,
\end{eqnarray*}
and so $\rho_\alpha\equiv\rho_\beta$ on $U_\alpha\cap U_\beta\cap S$ for all $\alpha$
and $\beta$ if and only if $\gt U$ is a splitting atlas.
\smallskip

(iii) Let $\gt U=\{(U_\alpha,z_\alpha)\}$ be a splitting atlas,
adapted to the first order lifting~$\rho_0$, and choose another
first order lifting~$\rho$. Lemma~\ref{th:uquattro}.(ii) implies that
$\rho=\rho_0-\iota\circ\phe$ for a suitable $\phe\in H^0\bigl(S,\Hom(\ca O_S,
\ca I_S/\ca I_S^2)\bigr)$, and it is easy to check that $\phe$ is a derivation. Therefore there
is a 0-cocycle
$\gt c=\{\gt c_\alpha\}\in H^0(\gt U_S,\ca N_S^*\otimes\ca T_S)$ such that
\[
\rho(g)=\rho_0(g)-(c_\alpha)^p_r\left[\frac{\de g}{\de z_\alpha^p} z_\alpha^r\right]_2
\]
for all $g\in\ca O_S|_{U_\alpha\cap S}$. Then defining new
coordinates as in (\ref{eq:ccuno}) we still get a splitting atlas,
easily seen adapted to~$\rho$. \end{proof}

\begin{remark}
Given a first order lifting $\rho\colon\ca O_S\to\ca O_{S(1)}$, let $\tilde\rho$,
$\tau^*$ and $\sigma^*$ be the morphisms associated to~$\rho$ by
Proposition~\ref{th:ucinque}, and let $\gt U=\{(U_\alpha,z_\alpha)\}$ be an atlas adapted to~$S$.
Then it is easy to check that the following assertions are equivalent:
\begin{itemize}
\item[(i)] $\gt U$ is adapted to~$\rho$;
\item[(ii)] for every $(U_\alpha,z_\alpha)\in\gt U$ with $U_\alpha\cap
S\ne\void$ and every $f\in\ca O_M|_{U_\alpha}$ one has
\[
\tilde\rho([f]_2)=\left[\frac{\de f}{\de z_\alpha^r}z_\alpha^r\right]_2;
\]
\item[(iii)] for every $(U_\alpha,z_\alpha)\in\gt U$ with $U_\alpha\cap
S\ne\void$ and every $r=1,\ldots,m$ one has
\[
\tau^*(\de_{r,\alpha})=\frac{\de}{\de z^r_\alpha};
\]
\item[(iv)] for every $(U_\alpha,z_\alpha)\in\gt U$ with $U_\alpha\cap
S\ne\void$ and every $f^j_\alpha\,\de/\de z^j_\alpha\in\ca T_{M,S}|_{U_\alpha\cap S}$ one
has
\[
\sigma^*\left(f^j_\alpha\frac{\de}{\de
z^j_\alpha}\right)=f^p_\alpha\frac{\de}{\de z^p_\alpha}.
\]
\end{itemize}
\end{remark}

\begin{remark}
\label{rem:linear}
It is also possible to prove (see, e.g.,~\cite{ABT3}) that a submanifold~$S$ splits into~$M$
if and only if its first infinitesimal neighbourhood~$S(1)$ in~$M$ is isomorphic to the
first infinitesimal neighbourhood~$S_N(1)$ of the zero section in~$N_S$.
\end{remark}

We end this section with a list of examples of splitting submanifolds.

\begin{example}
\label{ex:1}
A local holomorphic retract always splits in the ambient manifold (and thus it
is necessarily non-singular).
In particular, the zero section of a vector bundle always splits, as well as any slice
$S\times\{x\}$ in a product~$M=S\times X$ (with both $S$ and $X$ non-singular, of course).
\end{example}

\begin{example}
\label{ex:3}
If $S$ is a Stein submanifold of a complex manifold~$M$ (e.g., if $S$ is an open Riemann
surface), then $S$ splits into~$M$. Indeed, we have $H^1(S,\ca T_S\otimes\ca N_S^*)=(O)$ by
Cartan's Theorem~B, and the assertion follows from Proposition~\ref{th:unove}. In particular, if
$S$ is a \emph{singular} curve in $M$ then the non-singular part of~$S$ always splits in~$M$.
\end{example}

\begin{example}
\label{ex:4}
Let $\tilde M$ be the blow-up of a point in a complex manifold~$M$. Then the exceptional
divisor splits in~$\tilde M$: indeed, it is easy to check that the
atlas of~$\tilde M$ induced by the atlas of~$M$ is splitting.
\end{example}

\begin{example}
\label{ex:5}
A smooth closed irreducible subvariety of $\Pro^n$ splits into~$\Pro^n$ if and only if it is a linear
subspace (see \cite{VdV}, \cite{MP}, \cite{MR}).
\end{example}

\begin{example}
\label{ex:7}
Let $S$ be a non-singular, compact, irreducible curve of genus~$g$ on a surface~$M$.
If $S\cdot S<4-4g$ then $S$ splits into~$M$. In fact, the Serre duality for
Riemann surfaces implies that
\[
H^1\bigl(S,\Hom(\Omega_S,\ca N_S^*)\bigr)\cong H^0(S,\Omega_S\otimes\Omega_S\otimes\ca N_S),
\]
and the latter group vanishes because the line bundle $T^*S\otimes T^*S\otimes N_S$ has
negative degree by assumption. The bound $S\cdot S<4-4g$ is sharp: for instance, a
non-singular compact projective plane conic~$S$ has genus~$g=0$ and
self-intersection~$S\cdot S=4$, but it does not split in the projective plane (see
Example~\ref{ex:5}).
\end{example}

\section{Comfortably embedded submanifolds}
\label{sec:3}

In this section we introduce two other, more stringent, conditions on the embedding of $S$
into~$M$. From the definitions it will be clear that they are just the beginning of two infinite
lists of progressively more restrictive conditions; we shall however limit ourselves to present
only the properties we need in this paper, referring to~\cite{ABT3} for a more complete discussion.

We start with a natural generalization of splitting:

\begin{definition}
Let $S$ be a submanifold of a complex manifold~$M$. We say
that $S$ \emph{$2$-splits} into~$M$ if there exists a second order lifting~$\rho\colon\ca O_S\to\ca
O_{S(3)}$ or, in other words, if the exact sequence
\[
O\longrightarrow\ca I_S/\ca I_S^3\longhookrightarrow\ca O_{S(2)}
\mathop{\longrightarrow}\limits^{\theta_2}\ca O_S\longrightarrow O
\]
splits as sequence of sheaves of rings. Notice that a 2-splitting submanifold is necessarily
splitting.
\end{definition}

We have the following analogue of Proposition~\ref{th:udieci}:

\begin{proposition}
\label{th:pezzettino}
Let $S$ be an $m$-codimensional submanifold of a
complex manifold~$M$ of dimension~$n$. Then $S$ $2$-splits into~$M$ if and only if there is an
atlas $\gt U=\{(U_\alpha,z_\alpha)\}$ adapted to~$S$ such that
\begin{equation}
\frac{\de z_\beta^p}{\de z_\alpha^r}\in\ca I_S^2
\label{eq:ukappa}
\end{equation}
for all $r=1,\ldots,m$, $p=m+1,\ldots,n$ and indices $\alpha$, $\beta$ so that
$U_\alpha\cap U_\beta\cap S\ne\void$.
\end{proposition}

\begin{proof}
Let us first assume that we have an atlas $\gt U=\{(U_\alpha,z_\alpha)\}$ adapted to~$S$
and such that~(\ref{eq:ukappa}) holds. Let us then define $\rho_\alpha\colon\ca
O_S|_{U_\alpha}\to\ca O_{S(3)}|_{U_\alpha}$ by
\begin{equation}
\rho_\alpha([f]_1)=[f]_3-\left[\frac{\de f}{\de z_\alpha^r}z_\alpha^r\right]_3+
\frac{1}{2}\left[\frac{\de^2 f}{\de z_\alpha^{r_1}\de
z_\alpha^{r_2}}z_\alpha^{r_1}z_\alpha^{r_2}\right]_3
\label{eq:rhoalpha}
\end{equation}
for all $f\in\ca O(U_\alpha)$. It is easy to check that the right-hand side depends only
on~$[f]_1$, and that $\rho_\alpha$ is a ring morphism such that $\theta_2\circ\rho=\id$.
So to prove that $S$ 2-splits into~$M$ it suffices to show that $\rho_\alpha$ does not depend
on~$\alpha$. But indeed, since we have
\begin{equation}
[z_\beta^s]_3=\left[\frac{\de z^s_\beta}{\de z^r_\alpha}z_\alpha^r-
\frac{1}{2}\frac{\de^2 z^s_\beta}{\de z_\alpha^{r_1}\de z_\alpha^{r_2}}z_\alpha^{r_1}
z_\alpha^{r_2}\right]_3
\label{eq:pass}
\end{equation}
for all $s=1,\ldots,m$, we find
\begin{eqnarray}
\rho_\beta(\!\!\!&&[f]_1)-\rho_\alpha([f]_1)\nonumber\\
=&&\left[\frac{\de f}{\de z^r_\alpha}
z_\alpha^r\right]_3-\frac{1}{2}\left[\frac{\de^2 f}{\de z_\alpha^{r_1}\de
z_\alpha^{r_2}}z_\alpha^{r_1}z_\alpha^{r_2}\right]_3-
\left[\frac{\de f}{\de z^s_\beta} z_\beta^s\right]_3+\frac{1}{2}\left[\frac{\de^2 f}{\de
z_\beta^{s_1}\de z_\beta^{s_2}}z_\beta^{s_1}z_\beta^{s_2}\right]_3
\nonumber\\
=&&\left[\frac{\de f}{\de z_\beta^s}\left(\frac{\de z^s_\beta}{\de
z^r_\alpha}z_\alpha^r-z^s_\beta-\frac{1}{2}\frac{\de^2 z_\beta^s}{\de
z_\alpha^{r_1}\de z_\alpha^{r_2}}z_\alpha^{r_1}z_\alpha^{r_2}\right)\right]_3\nonumber\\
&&-\frac{1}{2}
\left[\frac{\de^2 f}{\de z_\beta^{s_1}\de z_\beta^{s_2}}\left(
\frac{\de z_\beta^{s_1}}{\de z_\alpha^{r_1}}\frac{\de z_\beta^{s_2}}{\de z_\alpha^{r_2}}
z_\alpha^{r_1}z_\alpha^{r_2}-z_\beta^{s_1}z_\beta^{s_2}\right)\right]_3\label{eq:uqbis}\\
&&+\left[\frac{\de f}{\de z^p_\beta}\frac{\de z^p_\beta}{\de z^r_\alpha}z_\alpha^r
\right]_3-\left[\frac{\de^2 f}{\de z^s_\beta\de z^p_\beta}\frac{\de z_\beta^s}{\de
z_\alpha^{r_1}}\frac{\de z_\beta^p}{\de z_\alpha^{r_2}} z_\alpha^{r_1}z_\alpha^{r_2}\right]_3
\nonumber\\
&&-\frac{1}{2}\left[\frac{\de^2 f}{\de z^{p_1}_\beta\de z^{p_2}_\beta}\frac{\de z_\beta^{p_1}}{\de
z_\alpha^{r_1}}\frac{\de z_\beta^{p_2}}{\de z_\alpha^{r_2}} z_\alpha^{r_1}z_\alpha^{r_2}\right]_3
\nonumber\\
=&&0\nonumber,
\end{eqnarray}
because of (\ref{eq:pass}) and (\ref{eq:ukappa}).

Conversely, let us assume that we have a second order lifting $\rho\colon\ca O_S\to\ca
O_{S(3)}$; we must build an atlas adapted to~$S$ such
that~(\ref{eq:ukappa}) holds.

If $\theta_{2,1}\colon\ca O_{S(2)}\to\ca O_{S(1)}$ is the canonical projection, then
$\rho_1=\theta_{2,1}\circ\rho$ is a first order lifting; let $\gt U=\{(U_\alpha,z_\alpha)\}$ be a
splitting atlas adapted to~$\rho_1$. Define then local second
order liftings~$\rho_\alpha$ as in~(\ref{eq:rhoalpha}), and set $\sigma_\alpha=\rho-\rho_\alpha$.
Now
$$
\theta_{2,1}\circ\sigma_\alpha=\rho_1-\theta_{2,1}\circ\rho_\alpha\equiv O,
$$
because the atlas is adapted to~$\rho_1$; therefore the image of~$\sigma_\alpha$ is
contained in~$\ca I_S^2/\ca I_S^3$, which is an~$\ca O_S$-module. Furthermore, $\sigma_\alpha$
is a derivation;
therefore we can find $(s_\alpha)^p_{r_1r_2}\in\ca O(U_\alpha\cap S)$, symmetric in
the lower indices, such that
\[
\sigma_\alpha=(s_\alpha)^p_{r_1r_2}[z_\alpha^{r_1}z_\alpha^{r_2}]_3\otimes\frac{\de}{\de
z_\alpha^p}.
\]
Now, since we are using a splitting atlas, the computations in (\ref{eq:uqbis}) yield
\begin{eqnarray}
\left[\frac{\de z^p_\beta}{\de z^r_\alpha}
z_\alpha^r\right]_3\otimes\frac{\de}{\de
z^p_\beta}&=&\rho_\beta-\rho_\alpha=\sigma_\alpha-\sigma_\beta\label{eq:uuuu}\\
&=&
\left[\left(\frac{\de z^p_\beta}{\de z^q_\alpha}(s_\alpha)^q_{r_1r_2}
-(s_\beta)^p_{s_1s_2}\frac{\de z_\beta^{s_1}}{\de z_\alpha^{r_1}}\frac{\de z_\beta^{s_1}}{\de
z_\alpha^{r_1}}\right)z_\alpha^{r_1}z_\alpha^{r_2}
\right]_3\otimes\frac{\de}{\de z^p_\beta}.\nonumber
\end{eqnarray}
Furthermore, always because we are using a splitting atlas, we can write
$z_\beta^p=\phi_{\beta\alpha}(z''_\alpha)+(h_{\beta\alpha})^p_{r_1r_2}z_\alpha^{r_1}z_\alpha^{r_2}$
with $(h_{\beta\alpha})^p_{r_1r_2}\in\ca O(U_\alpha\cap U_\beta)$ symmetric in the lower indices.
Putting this in (\ref{eq:uuuu}) we get
\[
2(h_{\beta\alpha})^p_{r_1r_2}-\frac{\de z^p_\beta}{\de z^q_\alpha}(s_\alpha)^q_{r_1r_2}
+(s_\beta)^p_{s_1s_2}\frac{\de z_\beta^{s_1}}{\de z_\alpha^{r_1}}\frac{\de z_\beta^{s_1}}{\de
z_\alpha^{r_1}}\in\ca I_S
\]
and hence
\begin{equation}
\left[\frac{\de z^p_\beta}{\de z^r_\alpha}-\frac{\de z^p_\beta}{\de z^q_\alpha}(s_\alpha)^q_{r_1r}
z_\alpha^{r_1}+(s_\beta)^p_{s_1s_2}z_\beta^{s_1}\frac{\de
z_\beta^{s_2}}{\de z_\alpha^r}
\right]_2=0.
\label{eq:uuud}
\end{equation}
Let us then consider the change of coordinates
\[
\begin{cases}\hat z_\alpha^r=z_\alpha^r,\cr \hat
z_\alpha^p=z_\alpha^p+\frac{1}{2}(s_\alpha)^p_{r_1r_2}(z''_\alpha)z_\alpha^{r_1}
z_\alpha^{r_2},\end{cases}
\]
defined in suitable open sets $\hat U_\alpha\subseteq U_\alpha$;
using (\ref{eq:uuud}) it is easy to check that $\{(\hat
U_\alpha,\hat z_\alpha)\}$ is the atlas we are looking for.
\end{proof}

\begin{definition}
An atlas adapted to $S$ satisfying (\ref{eq:ukappa}) will be said \emph{$2$-splitting.}
\end{definition}

\begin{remark}
\label{rem:cc2split}
There is a cohomological characterization of splitting submanifolds which are 2-splitting.
Indeed, assume that $S$ splits in~$M$, and let~$\gt U=\{(U_\alpha,z_\alpha)\}$ be a splitting
atlas. If we define local second order liftings~$\{\rho_\alpha\}$ as in~(\ref{eq:rhoalpha}), the
computations made in the proof of the previous proposition show that setting
$\rho_{\beta\alpha}=\rho_\beta-\rho_\alpha$ the cocycle $\{\rho_{\beta\alpha}\}$ defines a
cohomology class
\[
{\gt g}\in H^1(S,\ca I_S^2/\ca I_S^3\otimes\ca T_S),
\]
and it is easy to check that $S$ 2-splits if and only if this cohomology class vanishes.
\end{remark}

It turns out that for our aims it will be much more useful a different, though related, condition
on the embedding of~$S$ into~$M$.

Let $S$ be a splitting submanifold of codimension~$m$ of a
complex manifold~$M$, and let $\rho\colon\ca O_S\to\ca O_{S(1)}$ be a first order lifting.
The sheaf $\ca I_S/\ca I_S^3$ has a natural structure of $\ca O_{S(1)}$-module; by restriction of
scalars via~$\rho$, we get a structure of $\ca O_S$-module, and it is easy to check that with this
structure the sequence
\begin{equation}
O\longrightarrow\ca I_S^2/\ca I_S^3 \longhookrightarrow\ca I_S/\ca I_S^3
\mathop{\longrightarrow}\limits^{\theta_{2,1}}\ca I_S/\ca I_S^2 \longrightarrow O
\label{eq:ces}
\end{equation}
becomes an exact sequence of locally free $\ca O_S$-modules. In particular, it is clear that
if $\{(U_\alpha,z_\alpha)\}$ is an atlas adapted to~$\rho$ then
$\{[z_\alpha^r]_3,[z_\alpha^{r_1} z_\alpha^{r_2}]_3\}$ is a free set of local generators
of~$\ca I_S/\ca I_S^3$ over~$\ca O_S$.

\begin{remark}
\label{rem:comfemb}
The cohomology class ${\gt h}\in H^1\bigl(S,\ca I_S^2/\ca I_S^3\otimes\ca N_S\bigr)$ associated to
the sequence (\ref{eq:ces}) by the procedure described at the beginning of the previous section is
represented by the cocycle $\{{\gt h}_{\beta\alpha}\}$ given by
\[
{\gt h}_{\beta\alpha}=\frac{1}{2}\left[\frac{\de z_\beta^{s_1}}{\de z_\alpha^{r_1}}
\frac{\de z_\beta^{s_2}}{\de z_\alpha^{r_2}}\frac{\de^2 z_\alpha^r}{\de z^{s_1}_\beta
\de z^{s_2}_\beta}\right]_1\,[z_\alpha^{r_1}z_\alpha^{r_2}]_3\otimes\de_{r,\alpha},
\]
where $\{(U_\alpha,z_\alpha)\}$ is a splitting atlas associated to the first order lifting~$\rho$.
\end{remark}

We are thus led to the following

\begin{definition} Let $S$ be a (not necessarily closed) submanifold of a complex
manifold~$M$. We say that $S$ is \emph{comfortably embedded} in $M$ if there exists a first order
lifting $\rho\colon\ca O_S\to\ca O_{S(1)}$ such that the sequence (\ref{eq:ces})
splits as sequence of $\ca O_S$-modules. We shall sometimes say that $S$ is comfortably embedded
\emph{with respect to~$\rho$.}
\end{definition}

We can characterize comfortably embedded submanifolds using adapted atlases, recovering in
particular the original definition of comfortably embedded submanifolds introduced in~\cite{ABT}:

\begin{proposition}
\label{th:dtre}
Let $S$ be an $m$-codimensional submanifold of a
complex manifold~$M$ of dimension~$n$. Then $S$ is comfortably
embedded into~$M$ if and only if there exists an atlas $\gt
U=\{(U_\alpha,z_\alpha)\}$ adapted to~$S$ such that
\begin{equation}
\frac{\de z_\beta^p}{\de z_\alpha^r}\in\ca I_S\qquad\hbox{and}\qquad
\frac{\de^2 z_\beta^r}{\de z_\alpha^{s_1}\de z_\alpha^{s_2}}\in\ca I_S
\label{eq:ukappace}
\end{equation}
for all $r$,~$s_1$,~$s_2=1,\ldots,m$, $p=m+1,\ldots,n$ and indices $\alpha$, $\beta$ such that
$U_\alpha\cap U_\beta\cap S\ne\void$.
\end{proposition}

\begin{proof}
If we have an atlas satisfying (\ref{eq:ukappace}) then, by Proposition~\ref{th:udieci}, $S$ splits
into~$M$, and $\gt U$ is adapted to a first order lifting~$\rho$. Furthermore,
Proposition~\ref{th:uotto} and Remark~\ref{rem:comfemb} imply that $S$ is comfortably embedded
with respect to~$\rho$.

Conversely, assume that $S$ is comfortably embedded with respect to a first order lifting
$\rho\colon\ca O_S\to\ca O_{S(1)}$, let $\gt U=\{(U_\alpha,z_\alpha)\}$ be a splitting atlas
adapted to~$\rho$, and let $\nu\colon\ca I_S/\ca I_S^2\to\ca I_S/\ca I_S^3$ be a splitting morphism
for the sequence~(\ref{eq:ces}). For every index~$\alpha$ such that $U_\alpha\cap S\ne\void$
define~$\nu_\alpha\colon\ca I_S/\ca I_S^2|_{U_\alpha}\to\ca I_S/\ca I_S^3|_{U_\alpha}$ by setting
$\nu_\alpha([z_\alpha^r]_2)=[z_\alpha^r]_3$ and then extending by $\ca O_S$-linearity. Notice
that~(\ref{eq:defgrho}) and~(\ref{eq:pass}) yield
\begin{eqnarray}
\nu_\beta([z_\alpha^r]_2)&-&\nu_\alpha([z_\alpha^r]_2)\nonumber\\
&=&\rho\left(\left[\frac{\de
z_\alpha^r} {\de z_\beta^s}\right]_1\right)[z_\beta^s]_3-[z_\alpha^r]_3=-\frac{1}{2}
\left[\frac{\de^2 z^r_\alpha}{\de z_\beta^{s_1}\de z_\beta^{s_2}}z_\beta^{s_1}z_\beta^{s_2}
\right]_3.
\label{eq:uuut}
\end{eqnarray}
Now set $\sigma_\alpha=\nu-\nu_\alpha$; since $\theta_{2,1}\circ\sigma_\alpha\equiv O$, it follows
that $\IM\sigma_\alpha\subseteq\ca I_S^2/\ca I_S^3$. In particular, there are
$(c_\alpha)^r_{r_1r_2}\in\ca O(U_\alpha\cap S)$, symmetric in the lower indices, such that
\[
\sigma_\alpha([z^r_\alpha]_2)=(c_\alpha)^r_{r_1r_2}[z_\alpha^{r_1}z_\alpha^{r_2}]_3.
\]
But $\sigma_\alpha-\sigma_\beta=\nu_\beta-\nu_\alpha$; therefore (\ref{eq:uuut}) yields
\begin{equation}
(c_\beta)^s_{s_1s_2}\frac{\de z_\beta^{s_1}}{\de z_\alpha^{r_1}}\frac{\de z_\beta^{s_2}}{\de
z_\alpha^{r_2}}-\frac{\de z_\beta^s}{\de z_\alpha^r}(c_\alpha)^r_{r_1r_2}+
\frac{1}{2}\frac{\de^2 z^s_\beta}{\de z_\alpha^{r_1}\de z_\alpha^{r_2}}\in\ca I_S.
\label{eq:uuuq}
\end{equation}
We can finally define new coordinates $\hat z_\alpha$ by setting
\[
\begin{cases}\hat
z_\alpha^r=z_\alpha^r+(c_\alpha)^r_{r_1r_2}z_\alpha^{r_1}z_\alpha^{r_2},\cr
\hat z_\alpha^p=z_\alpha^p,\end{cases}
\]
on suitable open subsets $\hat U_\alpha\subseteq U_\alpha$. It is easy to check that $\hat{\gt
U}=\{(\hat U_\alpha,\hat z_\alpha)\}$ still is a splitting atlas adapted to~$\rho$. Moreover,
\[
\frac{\de^2\hat z^s_\beta}{\de\hat z_\alpha^{r_1}\hat z_\alpha^{r_2}}=
\frac{\de^2 z^s_\beta}{\de z_\alpha^{r_1}\de z_\alpha^{r_2}}+2\left((c_\beta)^s_{s_1s_2}
\frac{\de z_\beta^{s_1}}{\de z_\alpha^{r_1}}\frac{\de z_\beta^{s_2}}{\de
z_\alpha^{r_2}}-\frac{\de z_\beta^s}{\de z_\alpha^r}(c_\alpha)^r_{r_1r_2}
\right),
\]
and (\ref{eq:uuuq}) concludes the proof.
\end{proof}

\begin{definition}
An atlas satisfying (\ref{eq:ukappace}) will be said a {\sl comfortable
atlas.}
\end{definition}

We end this section with a last definition and some examples.

\begin{definition}
Let $S$ be a complex submanifold of a complex manifold~$M$.
We shall say that $S$ is \emph{$2$-linearizable} if it is 2-splitting and comfortably
embedded (with respect to the first order lifting induced by the 2-splitting).
\end{definition}

\begin{remark}
In \cite{ABT3} we prove that $S$ is 2-linearizable if and only if its second infinitesimal
neighbourhood~$S(2)$ in~$M$ is isomorphic to the second infinitesimal neighbourhood~$S_N(2)$
of the zero section in~$N_S$; compare with Remark~\ref{rem:linear}.
\end{remark}

\begin{example}
The zero section of a vector bundle is always 2-linearizable in the total space of the bundle.
\end{example}

\begin{example}
A local holomorphic retract is always 2-split in the
ambient manifold. Indeed, if $p\colon U\to S$ is a
local holomorphic retraction, then a second order lifting $\rho\colon\ca O_S\to\ca O_{S(2)}$ is
given by $\rho(f)=[f\circ p]_3$.
\end{example}

\begin{example}
Let $\tilde M$ be the blow-up of a submanifold~$X$ in a complex manifold~$M$. Then the
exceptional divisor~$E\subset\tilde M$ is 2-linearizable in~$\tilde M$.
\end{example}

\begin{example}
If $S$ is a Stein submanifold of a complex manifold~$M$ (e.g., if $S$ is an open Riemann
surface), then $S$ is 2-linearizable in~$M$.
Indeed, by Cartan's Theorem~B the first cohomology group of~$S$ with
coefficients in any coherent sheaf vanishes, and the assertion follows from
Proposition~\ref{th:uotto} and Remarks~\ref{rem:cc2split} and~\ref{rem:comfemb}. In particular, if
$S$ is a \emph{singular} curve in $M$ then the non-singular part of~$S$ is always comfortably
embedded in~$M$.
\end{example}

\begin{example}
Let $S$ be a non-singular, compact, irreducible curve of genus~$g$ in a surface~$M$. The
Serre duality for Riemann surfaces implies that
\[
H^1\bigl(S,\ca I_S^2/\ca I_S^3\otimes\ca T_S\bigr)\cong
H^0\bigl(S,\Omega_S\otimes\Omega_S\otimes\ca N_S^{\otimes 2}\bigr).
\]
Therefore if $2(S\cdot S)<4-4g$ then $H^1\bigl(S,\ca I_S^2/\ca I_S^3\otimes\ca T_S\bigr)
=(O)$. Analogously, we have
\[
H^1(S,\ca I_S^2/\ca I_S^3\otimes N_S)
\cong H^0(S,\Omega_S\otimes\ca N_S),
\]
and so $S\cdot S<2-2g$ implies $H^1(S,\ca I_S^2/\ca I_S^3\otimes N_S)=(O)$.
It follows that if $g\ge 1$ and $S\cdot S<4-4g$, or $g=0$ and $S\cdot S<2$, then
$S$ is 2-linearizable.
\end{example}

\section{Partial connections}
\label{sec:4}

As explained in the introduction, to get index theorems we need partial holomorphic connections.
Atiyah in~\cite{At} showed that a complex vector bundle admits a holomorphic connection if and only
if a particular exact sequence of locally free sheaves splits. In this section we shall adapt
Atiyah's construction to the case of partial holomorphic connections; in the next section we shall
describe a more concrete realization of Atiyah's exact sequence that will allow us to explicitely
construct splitting morphisms (the morphisms~$\psi$ of the introduction).

\begin{remark}
From now on, we shall denote the locally free sheaf of germs of holomorphic sections of a vector
bundle (e.g.,~$E$) by the corresponding calligraphic letter (e.g.,~$\ca E$).
\end{remark}

Let us start briefly recalling Atiyah's construction~\cite{At}. Let $E$ be a complex vector bundle
of rank $d$ over a complex manifold~$S$; we shall denote by $P_E$ the principal bundle
associated to $E$, with structure group~$GL(d)$. The group $GL(d)$ acts on the tangent vector bundle
$T_P$ of the total space of~$P_E$, and the quotient~$A_E=T_P/GL(d)$ can be identified with the
vector bundle on~$S$ of rank~$d^2-1$ composed by the fields of tangent
vectors to $P_E$ defined along one of its fibres and invariant under the action of $GL(d)$. Since
the action of~$GL(d)$ on~$P_E$ preserves the fibers of the canonical projection~$\pi_0\colon P_E\to
S$, the differential of~$\pi_0$ defines a vector bundle morphism, still denoted by $\pi_0$,
from~$A_E$ onto~$TS$. Atiyah has shown (\cite{At}, Theorem~1 and Proposition~9) that there is a
canonical exact sequence
\begin{equation}
O\longrightarrow\Hom(\ca E,\ca E) \longrightarrow \ca{A}_E\stackrel{\pi_0}{\longrightarrow} {\ca
T}_S\longrightarrow O,
\label{eq:extAE}
\end{equation}
of locally free $\ca O_S$-modules, where $\Hom(\ca E,\ca E)$ is canonically identified with the
sheaf of germs of holomorphic sections of the quotient, under the action of $GL(d)$, of the
sub-bundle of $T_P$ formed by vectors tangential to the fibres of~$P_E$.  Furthermore, this sequence
splits if and only if there is a holomorphic connection on~$E$ (\cite{At}, Theorem~2). See
also~\cite{GR}, where part of this theory is extended to subvarieties $S$
having normal crossing singularities.

\begin{remark}
\label{rem:extAE}
Atiyah (\cite{At}, p. 190 and~195) also computed the cohomology class
associated to the sequence~(\ref{eq:extAE}). In particular, if $E$ is
the normal bundle~$N_S$ of a submanifold~$S$ of a complex manifold~$M$ and
$\{(U_\alpha,z_\alpha)\}$ is an atlas adapted to~$S$, then the cohomology class is represented by
the cocycle $\{\gt g_{\alpha\beta}\}$ where
\begin{equation}
{\gt g}_{\alpha\beta}=\left.\frac{\de z^t_\beta}{\de z^s_\alpha}\frac{\de^2 z^r_\alpha}{\de
z^p_\beta\de z^t_\beta}\right|_S\,dz^p_\beta\otimes\omega^s_\alpha
\otimes\de_{r,\alpha}.
\label{eq:extAEb}
\end{equation}
\end{remark}

It is easy to adapt Atiyah's construction to the case of partial holomorphic connections.

\begin{definition}
Let $F$ be a sub-bundle of the tangent bundle~$TS$ of a complex manifold~$S$. A \emph{partial
holomorphic connection along~$F$} on a complex vector bundle~$E$ on~$S$ is a
$\C$-linear morphism $\nabla\colon\ca E\to\ca F^*\otimes\ca E$ such that
\[
\nabla(gs)=dg|_{\ca F}\otimes s+g\nabla s
\]
for all $g\in\ca O_S$ and $s\in\ca E$.
\end{definition}

If $F$ is a sub-bundle of~$TS$, we can consider the restriction to~$F$ of the
sequence~(\ref{eq:extAE})
\begin{equation}
O\longrightarrow\Hom(\ca E,\ca E) \longrightarrow \ca{A}_{E,F}\stackrel{\pi_0}{\longrightarrow} \ca
F\longrightarrow O,
\label{eq:extAF}
\end{equation}
where $\ca A_{E,F}=\pi_0^{-1}(\ca F)$. Then arguing as in \cite{At} it is easy to prove the
following

\begin{proposition}
\label{th:summ}
Let $F$ be a sub-bundle of the tangent bundle~$TS$ of a complex manifold~$S$, and let $E$ be a
complex vector bundle over~$S$. Then there is a partial holomorphic connection on~$E$ along~$F$ if
and only if the sequence~$(\ref{eq:extAF})$ splits, that is if and only if there is an $\ca
O_S$-morphism $\psi_0\colon\ca F\to\ca A_E$ such that~$\pi_0\circ\psi_0=\id$.
\end{proposition}

In the next section we shall give a more concrete realization of the sheaf $\ca A_E$
when $E$ is the normal bundle of a submanifold~$S$ into a manifold~$M$, allowing us to present an
alternative explicit description of the partial holomorphic connection given by a splitting of the
sequence~(\ref{eq:extAF}), and later on to build the morphisms~$\psi_0\colon\ca F\to\ca A_E$. But we
conclude this section with a few general definitions, useful to put in the right perspective what
we are going to do.

\begin{definition} Let $\ca E$ and $\ca F$ be locally free sheaves of $\ca O_S$-modules over
a complex manifold~$S$. Given a section~$X\in H^0(S,{\ca
T}_S\otimes\ca F^*)$, a \emph{holomorphic $X$-connection on~$\ca
E$\/} (also called a \emph{holomorphic action of $\ca F$ on $\ca E$
along~$X$\/}) is a $\C$-linear map $\tilde X\colon{\mathcal E}\to
{\mathcal F}^* \otimes {\mathcal E}$ such that
$$
\tilde{X}(gs)=X^*(dg)\otimes s+g\tilde{X}(s)
$$
for all $g\in\ca O_S$ and $s\in\ca E$, where $X^*\colon\Omega_S\to\ca F^*$ is the dual map of~$X$.
We shall often write
$\tilde X_v(s)$ instead of~$\tilde X(s)(v)$, for $s\in\ca E$ and~$v\in\ca F$, where as usual we have
identified~$\ca F^*\otimes\ca E$ with~$\Hom(\ca F,\ca E)$. Clearly,
if $\ca F$ is an $\ca O_S$-submodule of~$\ca T_S$ and $X$ is the
inclusion, then $\tilde X$ is just a partial holomorphic
connection on~$E$ along~$F$.
\end{definition}

\begin{remark}
Let $\ca E, \ca F$ be locally free $\ca O_S$-modules and $X\colon\ca F\to{\ca T}_S$ be an
$\ca O_S$-morphism. The \emph{sheaf of first jets
$J^1_X\ca E$ of $\ca E$ along~$X$} is the sheaf of abelian groups $(\ca
F^*\otimes \ca E)\oplus \ca E$, with the
structure of $\ca O_S$-module given  as follows: for $f\in \ca O_S$ and
$(\omega\otimes e)\oplus e' \in J^1_X \ca E$ we define
\[
f
(\omega\otimes e)\oplus e'=(X^*(df)\otimes e'+\omega \otimes fe)\oplus
fe'.
\]
The natural projection $J^1_X \ca E\to \ca E$ given by
$(\omega\otimes e)\oplus e'\to e'$ is a surjective $\ca
O_S$-morphism whose kernel is $\ca F^*\otimes \ca E$.
Thus we
obtain the exact sequence
\begin{equation}
O\longrightarrow \ca F^*\otimes \ca E\longrightarrow J^1_X \ca E\longrightarrow \ca E\to O.
\label{eq:jet}
\end{equation}
Notice that $J^1_X \ca E$ is locally
$\ca O_S$-free and the sequence (\ref{eq:jet})
 is functorial on $\ca E$. It is easy to see that this sequence splits if and only if
there is a holomorphic $X$-connection on~$E$. If we denote by $c_X(\ca E)\in H^1\bigl(S,{\ca
F}^*\otimes\Hom(\ca E,\ca E)\bigr)$ the class associated to the sequence~(\ref{eq:jet}),
and by $\hat c(\ca E)\in H^1\bigl(S,\Omega_S\otimes\Hom(\ca E,\ca E)\bigr)$ the class associated to
the same sequence when ${\ca F}=\ca T_S$ and $X$ is the identity, then it is not difficult to see
that
\[
c_X(\ca E)=(X^*\otimes \id)_*\hat c({\mathcal E}).
\]
Furthermore, Atiyah (\cite{At}, Theorem~5) has shown that~$\hat c(\ca E)$ is the opposite of the
cohomology class associated to the sequence~(\ref{eq:extAE}).
\end{remark}

We shall need a notion of flatness for a holomorphic $X$-connection. To state it in full
generality, we need a new definition and a lemma.

\begin{definition}
Let $\ca F$ be a sheaf of $\ca O_S$-modules over a complex manifold~$S$, equipped
with a $\ca O_S$-morphism $X\colon\ca F\to\ca T_S$. We say that $\ca F$ is a \emph{Lie
algebroid} of \emph{anchor~$X$} if there is a $\C$-bilinear map $\{\cdot\,,\cdot\}\colon\ca
F\oplus\ca F\to\ca F$ such that
\begin{itemize}
\item[(a)] $\{v,u\}=-\{u,v\}$;
\item[(b)] $\{u,\{v,w\}\}+\{v,\{w,u\}\}+\{w,\{u,v\}\}=O$;
\item[(c)] $\{g\cdot u,v\}=g\cdot\{u,v\}-X(v)(g)\cdot u$ for all $g\in\ca O_S$ and
$u$,~$v\in\ca F$.
\end{itemize}
\end{definition}

\begin{example} Assume that $X\colon\ca F\to\ca T_S$ is injective, and that
$X(\ca F)$ is an \emph{involutive} subsheaf of~$\ca T_S$ (we recall
that a subsheaf $\ca F$ of $\ca T_S$ is involutive if it is locally
$\ca O_S$-free and, for each $x\in S$, the fiber $\ca F_x$ is closed
under the bracket operation for vector fields). Then we can easily
provide $\ca F$ with a Lie algebroid structure of anchor~$X$ by setting
$$
\{u,v\}=X^{-1}([X(u),X(v)]).
$$
In particular in this case we have $X(\{u,v\})=[X(u),X(v)]$.
\end{example}

We refer to \cite{M} and references therein for the general theory of Lie algebroids; here we shall
use the definition only.

Let $\ca F$ be a Lie algebroid of anchor~$X$ over a complex
manifold~$S$, and assume we also have a holomorphic
$X$-connection $\tilde X\colon\ca E\to\ca F^*\otimes\ca E$ over a locally free $\ca
O_S$-module~$\ca E$. Then it is easy to see that setting
\begin{equation}
R_{uv}s=\tilde X_u\circ\tilde X_v(s)-\tilde X_v\circ\tilde X_u(s)-\tilde X_{\{u,v\}}(s)
\label{eq:nuovouno}
\end{equation}
we define a $\C$-linear map $R\colon\ca E\to\bigwedge^2\ca F^*\otimes\ca E$ such that
\begin{equation}
R_{uv}(g\cdot s)=g\cdot R_{uv}(s)+\bigl([X(u),X(v)]-X(\{u,v\})\bigr)(g)\cdot s
\label{eq:nuovodue}
\end{equation}
for all $g\in\ca O_S$, $u$,~$v\in\ca F$ and $s\in\ca E$.

\begin{definition} Let $\ca F$ be a Lie algebroid of anchor~$X$ over a complex
manifold~$S$, and $\tilde X\colon\ca E\to\ca F^*\otimes\ca E$ a
holomorphic $X$-connection over a locally free $\ca O_S$-module~$\ca E$. The \emph{curvature}
of~$\tilde X$ (with respect to the given Lie algebroid structure) is the $\C$-linear map
$R\colon\ca E\to\bigwedge^2\ca F^*\otimes\ca E$ defined in (\ref{eq:nuovouno}). We shall say that
$\tilde X$ is \emph{flat} if $R\equiv O$.
\end{definition}

\begin{remark}
Note that if $\tilde X$ is flat then from (\ref{eq:nuovodue}) it
follows that $[X(u),X(v)]=X(\{u,v\})$ for all $u,v\in \ca F$.
\end{remark}

\begin{example} Let $X\colon\ca F\to \ca T_S$ be the inclusion of an involutive locally free $\ca
O_S$-submodule $\ca F$ of~$\ca T_S$ with locally free quotient
$\ca Q=\ca T_S/\ca F$ and consider the  associated exact
sequence
\[
O\longrightarrow\ca F\stackrel{X}{\longhookrightarrow}\ca T_S \stackrel{\phi}{\longrightarrow}\ca
Q\longrightarrow O.
\]
Then we can define a partial holomorphic connection $\tilde X\colon\ca Q\to\ca
F^*\otimes\ca Q$ on~$\ca Q$ along~$\ca F$ by setting
$$
\tilde X_u(q)=\phi([X(u),\tilde q]),
$$
where $\tilde q$ is any local section of~$\ca T_S$ such that
$\phi(\tilde q)=q$. Putting on $\ca F$ the natural Lie algebroid structure of anchor~$X$ given
by the bracket of vector fields, it is easy to check that $\tilde X$ is a
flat partial holomorphic connection. The existence of this
natural flat partial holomorphic connection is one of the main ingredients in the
proof of Baum-Bott index theorem for singular holomorphic
foliations; see~\cite{BB}, \cite{CL}, and \cite{S2}, Chapter~I\negthinspace
I\negthinspace I.
\end{example}

\section{The universal partial connection on the normal sheaf}
\label{sec:5}
In this section we shall describe a concrete incarnation of the sheaf~$\ca A_E$ when $E$ is the
normal bundle of a submanifold~$S$ of a complex manifold~$M$.

\begin{definition}
Let $S$ be a (not necessarily closed) complex submanifold of a manifold~$M$. Set $\ca T_{M,S(1)}=\ca
T_M\otimes_{\ca O_M}\ca O_{S(1)}$; with a slight abuse of notation we shall denote
by~$\theta_1$ the
$\ca O_M$-morphism~$\id\otimes\theta_1\colon\ca T_{M,S(1)}\to\ca T_{M,S}$.
Let $\ca
T^S_{M,S(1)}$ be the $\ca O_M$-submodule of~$\ca T_{M,S(1)}$ given by
\[
\ca T^S_{M,S(1)}=\ker(p_2\circ\theta_1)\subset\ca T_{M,S(1)},
\]
where $p_2\colon\ca T_{M,S}\to\ca N_S$ is the natural projection.
Notice that $\theta_1(\ca T^S_{M,S(1)})=\ca T_S$.
\end{definition}

In local adapted coordinates, an element $v=[a^j]_2\frac{\de}{\de z^j}\in\ca T_{M,S(1)}$
belongs to~$\ca T^S_{M,S(1)}$ if and only if $[a^r]_1=0$ for $r=1,\ldots,m$. In other
words, $v\in\ca T_{M,S(1)}$ belongs to~$\ca T^S_{M,S(1)}$ if and only if when restricted
to~$S$ it is tangent to it.

In general, $\ca T^S_{M,S(1)}$ is not an $\ca O_S$-module (it is if $S$ splits into~$M$, but we
do not want to assume this yet). However, we can almost define on it a Lie algebroid structure of
anchor~$\theta_1$.

If $v\in\ca T_{M,S}$ and $f\in\ca O_M$, then $v(f)$ is a
well-defined element of~$\ca O_S$. Analogously, if $v\in\ca T_{M,S(1)}$ and $f\in\ca
O_M$, then $v(f)$ is a well-defined element of~$\ca O_{S(1)}$; on the other hand,
if~$g\in\ca O_{S(1)}$ then $v(g)$ is well defined as an element
of~$\ca O_S$, but not of~$\ca O_{S(1)}$. This means that we can define a bracket
operation $[\cdot\,,\cdot]\colon\ca T_{M,S(1)}\oplus\ca
T_{M,S(1)}\to\ca T_{M,S}$ by setting
\[
[u,v](f)=u\bigl(v(f)\bigr)-v\bigl(u(f)\bigr)\in\ca O_S
\]
for all $f\in\ca O_M$.
In particular, for every $g\in\ca O_{S(1)}$ and $u$,~$v\in\ca T_{M,S(1)}$ we have
\begin{equation}
[g u,v]=\theta_1(g)[u,v]-v(g)\cdot\theta_1(u).
\label{eq:nuovozer}
\end{equation}

\begin{remark}
In general, if $u$,~$v$,~$w\in\ca T_{M,S(1)}$ then $[v,w]\in\ca T_{M,S}$, and so
$[u,[v,w]]$ is {\it not} defined. But we shall see exceptions to this rule.
\end{remark}

\begin{lemma}
\label{th:nuovodue}
Let $S$ be a complex submanifold of a manifold~$M$. Then
\begin{itemize}
\item[\textup{(i)}]
every $v\in {\mathcal T}^S_{M,S(1)}$ induces a derivation $g\mapsto
v(g)$ of~$\ca O_{S(1)}$;
\item[\textup{(ii)}]
there exists a natural $\C$-bilinear map~$\{\cdot\,,\cdot\}\colon\ca
T^S_{M,S(1)}\oplus \ca T^S_{M,S(1)}\to\ca T^S_{M,S(1)}$ such that
\begin{itemize}
\item[\textup{(a)}] $\{v,u\}=-\{u,v\}$,
\item[\textup{(b)}] $\{u,\{v,w\}\}+\{v,\{w,u\}\}+\{w,\{u,v\}\}=O$,
\item[\textup{(c)}] $\{gu,v\}=g\{u,v\}-v(g)\cdot u$ for all $g\in\ca O_{S(1)}$, and
\item[\textup{(d)}] $\theta_1 \{u,v\}=[\theta_1(u),\theta_1(v)]=[u,v]$.
\end{itemize}
\end{itemize}
\end{lemma}

\begin{proof}
(i) The relevant fact here is that if we have an adapted chart $(U,z)$
and germs $[f]_2\in\ca O_{S(1)}$ and~$[h]_2\in\ca I_S/\ca I_S^2$, then we get well-defined
elements of~$\ca O_{S(1)}$ by setting
\[
\frac{\de[f]_2}{\de z^p}=\left[\frac{\de f}{\de z^p}\right]_2\qquad\hbox{and}\qquad
[h]_2\frac{\de[f]_2}{\de z^r}=\left[h\frac{\de f}{\de z^r}\right]_2
\]
for $p=m+1,\ldots,n$ and $r=1,\ldots,m$. This implies that if $[f]_2\in\ca O_{S(1)}$ and
$v=[a^j]_2\frac{\de}{\de z^j}\in\ca T^S_{M,S(1)}$ then
\[
v([f]_2)=[a^r]_2\frac{\de [f]_2}{\de z^r}+[a^p]_2\frac{\de [f]_2}{\de z^p}
\]
is a well-defined (and independent of the local coordinates) element of~$\ca O_{S(1)}$,
and not just of~$\ca O_S$, and in this way we clearly get a derivation of~$\ca O_{S(1)}$.
\smallskip

(ii) We define $\{\cdot\,,\cdot\}$ by setting
\[
\{u,v\}(g)=u\bigl(v(g)\bigr)-v\bigl(u(g)\bigr),
\]
for all $g\in\ca O_{S(1)}$. It is easy to check (working for
instance in local coordinates adapted to~$S$) that $\{u,v\}\in\ca
T^S_{M,S(1)}$, and that properties (a)--(d) are satisfied.
\end{proof}

It turns out that a quotient of~$\ca T^S_{M,S(1)}$ has
a natural $\ca O_S$-module structure, and inherites
the Lie algebroid structure of anchor~$\theta_1$.  To prove this, we
need the following

\begin{lemma}
\label{th:nuovotre}
Let $S$ be an $m$-codimensional complex submanifold of a manifold~$M$ of dimension~$n$. Then:
\begin{itemize}
\item[\textup{(i)}] $u\in\ca T_{M,S(1)}$ is such that $p_2([u,s])=O$ for all $s\in\ca T_{M,S(1)}$ if
and only if $u\in\ca I_S\cdot\ca T^S_{M,S(1)}$;
\item[\textup{(ii)}] if $u\in\ca I_S\cdot\ca T^S_{M,S(1)}$ and $v\in\ca T^S_{M,S(1)}$ then
$\{u,v\}\in\ca I_S\cdot\ca T^S_{M,S(1)}$;
\item[\textup{(iii)}]the quotient sheaf
\[
\ca A=\ca T^S_{M,S(1)}/\ca I_S\cdot\ca T^S_{M,S(1)}
\]
admits a natural structure of $\ca O_S$-locally free sheaf such that $\theta_1\colon\ca{A}\to\ca
T_S$ is an $\ca O_S$-morphism.
\end{itemize}
\end{lemma}

\begin{proof} (i) Let us work in local coordinates adapted to~$S$. Writing
$u=[a^j]_2\frac{\de}{\de z^j}$, $s=[b^h]_2\frac{\de}{\de z^h}\in\ca T_{M,S(1)}$ we have
\[
p_2([u,s])=\left[a^j\frac{\de b^r}{\de z^j}-b^j\frac{\de a^r}{\de z^j}\right]_1\frac{\de}{\de
z^r}.
\]
Now $u\in\ca I_S\cdot\ca T^S_{M,S(1)}$ if and only if $a^r\in\ca I_S^2$ and $a^p\in\ca
I_S$ for $r=1,\ldots,m$ and $p=m+1,\ldots,n$; in particular it is clear that $u\in\ca
I_S\cdot\ca T^S_{M,S(1)}$ implies $p_2([u,s])=O$ for all~$s\in\ca T_{M,S(1)}$.

Conversely, assume that $u=[a^j]_2\frac{\de}{\de z^j}\in\ca T_{M,S(1)}$ is such that
$p_2([u,s])=O$ for all $s\in\ca T_{M,S(1)}$. From $p_2([u,\de/\de z^p])=O$ for all
$p=m+1,\ldots, n$ we get that $[a^r]_1$ is a constant~$\alpha^r\in\C$ for $r=1,\ldots,m$.
But then from $p_2([u,\de/\de z^s])=O$ for all
$s=1,\ldots, m$ we get $[a^r]_2=\alpha^r$ for $r=1,\ldots,m$. Now from
$p_2([u,[z^{s_0}]_2\de/\de z^{s_0}])=O$ for all
$s_0=1,\ldots, m$ (no sum on~$s_0$ here) we get $\alpha_r=0$ for $r=1,\ldots,m$. Finally,
from $p_2([u,[z^{p_0}]_2\de/\de z^1])=O$ for all
$p_0=m+1,\ldots, n$ we get $[a^p]_1=0$ for all
$p=m+1,\ldots, n$, and so~$u\in\ca I_S\cdot\ca T^S_{M,S(1)}$, as claimed.
\smallskip

(ii) Working again in local coordinates adapted to~$S$, if $u=[a^j]_2\frac{\de}{\de
z^j}\in\ca I_S\cdot\ca T^S_{M,S(1)}$ and $v=[b^h]_2\frac{\de}{\de z^h}\in\ca T^S_{M,S(1)}$
then we have
\[
\{u,v\}=\left[a^j\frac{\de b^r}{\de z^j}-b^j\frac{\de a^r}{\de
z^j}\right]_2\frac{\de}{\de z^r} +\left[a^j\frac{\de b^p}{\de z^j}-b^j\frac{\de a^p}{\de
z^j}\right]_2\frac{\de}{\de z^p}\in\ca I_S\cdot\ca T^S_{M,S(1)}
\]
because $a^r\in\ca I_S^2$ and $a^p$,~$b^r\in\ca I_S$ for all $r=1,\ldots,m$ and
$p=m+1,\ldots,n$.
\smallskip

(iii) The sheaf $\ca T^S_{M,S(1)}$ is an $\ca O_{S(1)}$-submodule of~$\ca T_{M,S(1)}$
such that (by definition) $g\cdot v\in\ca I_S\cdot\ca T^S_{M,S(1)}$ for all $g\in\ca
I_S/\ca I_S^2$ and $v\in\ca T^S_{M,S(1)}$. Therefore the $\ca O_{S(1)}$-module structure
induces a natural $\ca O_S$-module structure on~$\ca{A}$. It is easy to check that,
in terms of local coordinates adapted to $S$, the sheaf ${\ca A}$ is a locally free
$\ca O_S$-module freely generated by $\pi(\frac{\de}{\de z^p})$ and $\pi([z^s]_2\frac{\de}{\de
z^r})$ (with $p=m+1, \ldots n$, and
$r$,~$s = 1,\ldots m$), where $\pi\colon\ca T^S_{M,S(1)}\to\ca A$ is the quotient map.

Finally, since $\ca I_S\cdot\ca T^S_{M,S(1)}\subset\ker\theta_1$, the morphism $\theta_1$ defines a
map, still denoted by the same symbol, $\theta_1\colon\ca{A}\to\ca T_S$, and it
is clear that $\theta_1$ is an $\ca O_S$-morphism for the structure we just defined.

\end{proof}

\begin{definition} Let $S$ be a complex submanifold of a complex manifold~$M$.
The \emph{Atiyah sheaf} of~$S$ in~$M$ is the locally free $\ca O_S$-module
\[
\ca{A}=\ca T^S_{M,S(1)}/\ca I_S\cdot\ca T^S_{M,S(1)}.
\]
\end{definition}

The sheaf $\ca{A}$ appears also in \cite{PS}, where it is denoted by~${\ca
N}^1_{M,C}$, and in~\cite{GR}, where it is denoted by~$T^1_M\langle Y\rangle\otimes\ca O_Y$. We now
show that $\ca A$ is isomorphic to the sheaf~$\ca A_{N_S}$ described in the previous section.

\begin{theorem}
\label{th:estensioneA}
Let $S$ be an $m$-codimensional complex submanifold of a manifold~$M$ of dimension~$n$. Then
there is a natural exact sequence of locally free $\ca O_S$-modules
\begin{equation}
 O\longrightarrow\Hom({\ca N}_S,\ca N_S) \longrightarrow \ca{A}\stackrel{\theta_1}{\longrightarrow}
{\ca T}_S\longrightarrow O
\label{eq:extA}
\end{equation}
which is isomorphic to the sequence~$(\ref{eq:extAE})$ with $E=N_S$. In particular, the sheaf~$\ca
A$ only depends on the normal bundle $N_S$ of the embedding of $S$ into~$M$, and $N_S$ admits
a holomorphic connection if and only if the sequence~$(\ref{eq:extA})$ splits.
\end{theorem}

\begin{proof}
As usual, we work in local coordinates~$\{(U_\alpha,z_\alpha)\}$ adapted to~$S$.
The kernel of $\theta_1$ is freely generated by the images under the canonical
projection~$\pi\colon\ca T^S_{M,S(1)}\to\ca A$ of $[z^s_\alpha]_2\frac{\de}{\de z^r_\alpha}$ for
$r$,~$s=1,\ldots,m$.  Now we have
\[
\pi\left([z^s_\alpha]_2\frac{\de}{\de z^r_\alpha}\right)=\left.\frac{\de
z^s_\alpha}{\de z^{s_1}_\beta}\right|_{S}
\left.\frac{\de z^{r_1}_\alpha}{ \de z^r_\beta}\right|_{S}
\pi\left([z^{s_1}_\beta]_2 \frac{\de}{\de z^{r_1}_\beta}\right),
\]
and hence the kernel of
$\theta_1$ is naturally isomorphic to $\Hom(\ca N_S,\ca N_S)$.

To prove that (\ref{eq:extA}) is isomorphic to (\ref{eq:extAE}), by Proposition~\ref{th:uotto}
it suffices to prove that the cohomology class associated to both sequences is the same.

Define local splittings~$\sigma_\alpha$ of (\ref{eq:extA}) by setting $\sigma_\alpha(\frac{\de}{\de
z^p_\alpha}) =\pi\left(\frac{\de}{\de z^p_\alpha}\right)$ and then extending by $\ca
O_S$-linearity.
The class of the sequence (\ref{eq:extA}) is then represented by the cocycle given by
\begin{eqnarray*}
(\sigma_\beta-\sigma_\alpha)\left(\frac{\de}{\de z^p_\beta}\right)
&=&\pi\left(\frac{\de}{\de z^p_\beta}\right)
-\left.\frac{\de z^q_\alpha}{\de z^p_\beta}\right|_{S}
\pi\left(\frac{\de}{\de z^q_\alpha}\right)\\
&=& \pi\left(\left[\frac{\de z_\alpha^r}{\de z_\beta^p}\right]_2
\frac{\de}{\de z_\alpha^r}\right)=
\pi\left(\left[\frac{\de^2 z_\alpha^r}{\de z_\beta^p\de z_\beta^t}z_\beta^t\right]_2
\frac{\de}{\de z_\alpha^r}\right)\\
&=&\left.\frac{\de z_\beta^t}{\de z_\alpha^s}\frac{\de^2 z_\alpha^r}{\de z_\beta^p\de z_\beta^t}
\right|_S
\pi\left([z_\alpha^s]_2\frac{\de}{\de z_\alpha^r}
\right),
\end{eqnarray*}
and our claim follows from~(\ref{eq:extAEb}). The last assertion now
is an immediate consequence of \cite{At}, Theorem~2.
\end{proof}

The main advantage of our~$\ca A$ over Atiyah's~$\ca A_{N_S}$ is that we have
an explicit way of going from a splitting of the sequence~(\ref{eq:extA}) to a partial holomorphic
connection on~$N_S$. Indeed, $\ca A$ comes equipped with both a natural structure of Lie algebroid
of anchor~$\theta_1$ and a holomorphic $\theta_1$-connection on~$N_S$:

\begin{proposition}
\label{th:nuovoqua}
Let $S$ be a complex submanifold of a complex
manifold~$M$. Then:
\begin{itemize}
\item[\textup{(i)}] the Atiyah sheaf $\ca{A}$ has a natural structure~$\{\cdot\,,\cdot\}$
of Lie algebroid of anchor~$\theta_1$ such that
\begin{equation}
\theta_1\{q_1,q_2\}=[\theta_1(q_1),\theta_1(q_2)]
\label{eq:nuovotre}
\end{equation}
for all~$q_1$,~$q_2\in\ca{A}$;
\item[\textup{(ii)}] there is a natural holomorphic $\theta_1$-connection $\tilde
X\colon\ca N_S\to\ca{A}^*\otimes\ca N_S$ on~$\ca N_S$ given by
\[
\tilde X_q(s)=p_2([v,\tilde s])
\]
for all $q\in\ca{A}$ and $s\in\ca N_S$, where $v\in\ca T^S_{M,S(1)}$ and $\tilde s\in\ca
T_{M,S(1)}$ are such that $\pi(v)=q$ and $p_2\circ\theta_1(\tilde s)=s$;
\item[\textup{(iii)}] this holomorphic $\theta_1$-connection $\tilde X$ is flat.
\end{itemize}
\end{proposition}

\begin{proof}
(i) Lemmas~\ref{th:nuovodue} and~\ref{th:nuovotre}.(ii) imply that setting
\[
\{q_1,q_2\}=\pi(\{v_1,v_2\})
\]
for all $q_1$,~$q_2\in\ca{A}$, where $v_j\in\ca T^S_{M,S(1)}$ is such that
$q_j=\pi(v_j)$, we get a well-defined $\theta_1$-Lie algebroid structure
satisying (\ref{eq:nuovotre}).
\smallskip

(ii) Lemma~\ref{th:nuovotre}.(i) implies that if $\pi(v)=\pi(v')$ then $p_2([v,\tilde
s])=p_2([v',\tilde s])$ for all $\tilde s\in\ca T_{M,S(1)}$. Analogously,
Lemma~\ref{th:nuovodue}.(ii).(d) implies that if $p_2\circ\theta_1(\tilde
s)=p_2\circ\theta_1(\tilde s')$ then $p_2([v,\tilde
s])=p_2([v,\tilde s'])$ for all $v\in\ca T^S_{M,S(1)}$; therefore $\tilde X_q(s)$ is a
well-defined element of~$\ca N_S$.

Now, (\ref{eq:nuovozer}) yields
\[
\tilde X_{[f]_1\cdot q}(s)=p_2\bigl(\bigl[[f]_2v,\tilde s\bigr]\bigr)=
p_2\bigl([f]_1[v,\tilde s]-\tilde s([f]_2)\cdot\theta_1(v)\bigr)=[f]_1\cdot
\tilde X_q(v)
\]
because $v\in\ker(p_2\circ\theta_1)$, and so $\tilde X(s)\in\ca{A}^*\otimes\ca N_S$,
as claimed. Finally,
\begin{eqnarray*}
\tilde X_q([f]_1\cdot s)&=&p_2\bigl(\bigl[v,[f]_2\tilde s\bigr]\bigr)=
p_2\bigl([f]_1[v,\tilde s]+v([f]_2)\cdot\theta_1(\tilde s)\bigr)\\
&=&[f]_1\cdot\tilde X_q(v)
+\theta_1(q)([f]_1)\cdot s,
\end{eqnarray*}
because
$$
\bigl[v([f]_2)\bigr]_1=\theta_1(v)([f]_1)=\theta_1\bigl(\pi(v)\bigr)([f]_1)
$$
for all $[f]_2\in\ca O_{S(1)}$ and $v\in\ca T^S_{M,S(1)}$, and thus $\tilde X$ is a
holomorphic $\theta_1$-connection.
\smallskip

(iii) We must prove that
\begin{equation}
\bigl[v_1,\widetilde{[v_2,\tilde s]}\bigr]+\bigl[v_2,\widetilde{[\tilde s,v_1]}\bigr]
+[\tilde s,\{v_1,v_2\}]\in\ca T_S
\label{eq:nuovoqua}
\end{equation}
for all $v_1$, $v_2\in\ca T^S_{M,S(1)}$ and $\tilde s\in\ca
T_{M,S(1)}$, where $\widetilde{[v_j,\tilde s]}$ is any element
of~$\ca T_{M,S(1)}$ such that its~$\theta_1$-image is equal
to~$[v_j,\tilde s]\in\ca T_{M,S}$. But using local coordinates
adapted to~$S$ it is easy to see that (\ref{eq:nuovoqua}) is a
consequence of the usual Jacobi rule for brackets of vector fields.
\end{proof}

\begin{definition} Let $S$ be a complex submanifold of a complex manifold~$M$. The
holomorphic $\theta_1$-connection $\tilde X\colon\ca N_S\to\ca{A}^*\otimes\ca N_S$ on~$\ca N_S$
defined in Proposition~\ref{th:nuovoqua}.(ii) is called the \emph{universal holomorphic connection}
on~$\ca N_S$.
\end{definition}

We can now summarize what we have done up to now in the following

\begin{theorem}
\label{th:sommario}
Let $S$ be a submanifold of a complex manifold~$M$, and $F$ a sub-bundle of the tangent
bundle~$TS$. Then:
\begin{itemize}
\item[\textup{(i)}] if $\psi\colon\ca F\to\ca A$ is an $\ca O_S$-morphism such that
$\theta_1\circ\psi=\id$ then the map $\tilde X^\psi\colon\ca N_S\to\ca F^*\otimes\ca N_S$ given by
\[
\tilde X^\psi_v(s)=\tilde X_{\psi(v)}(s)
\]
for all $v\in\ca F$ and $s\in\ca N_S$, where $\tilde X$ is the universal holomorphic connection
on~$N_S$, is a partial holomorphic connection on~$N_S$ along~$F$;
\item[\textup{(ii)}] there exists a partial holomorphic connection on~$N_S$ along~$F$ if and only if
there exists an $\ca O_S$-morphism $\psi\colon\ca F\to\ca A$ such that $\theta_1\circ\psi=\id$;
\item[\textup{(iii)}] if $F$ is involutive, then the partial holomorphic connection $\tilde X^\psi$
is flat if and only if $\psi\colon\ca F\to\ca A$ is a Lie algebroid morphism.
\end{itemize}
\end{theorem}

\begin{proof}
(i) The only not completely trivial property is Leibniz's rule. But indeed
\[
\tilde X^\psi_v(g\cdot s)=\tilde X_{\psi(v)}(g\cdot s)=g\cdot\tilde X_{\psi(v)}(s)+
\theta_1\bigl(\psi(v)\bigr)(g)\cdot s=g\cdot\tilde X^\psi_v(s)+v(g)\cdot s
\]
for all $g\in\ca O_S$, and we are done.
\smallskip

(ii) In one direction is (i). Conversely, assume that we have a partial holomorphic connection
on~$N_S$ along~$F$. Then Proposition~\ref{th:summ} yields an $\ca O_S$-morphism~$\psi_0$ from~$\ca
F$ to~$\ca A_{N_S}$ such that $\pi_0\circ\psi_0=\id$, and hence Theorem~\ref{th:estensioneA} yields
an
$\ca O_S$-morphism $\psi\colon\ca F\to\ca A$ such that $\theta_1\circ\psi=\id$.
\smallskip

(iii) If we denote by~$R^\psi$
the curvature of~$\tilde X^\psi$, recalling that the universal
holomorphic connection~$\tilde X$ is flat, we get
\[
R^\psi_{uv}=\tilde
X_{\{\psi(u),\psi(v)\}-\psi[u,v]}
\]
for all $u$,~$v\in\ca F$. One direction is then clear; conversely, assume that $R^\psi\equiv O$.
Now, Proposition~\ref{th:nuovoqua}.(i) implies that
$\{\psi(u),\psi(v)\}-\psi[u,v]\in\ker\theta_1\subset\ca{A}$ for all~$u$,~$v\in\ca F$;
therefore it suffices to prove that if $q\in\ker\theta_1\subset\ca{A}$ is such
that~$\tilde X_q\equiv O$ then $q=O$. But indeed, if $q\in\ker\theta_1$ then $q=\pi(v)$
with~$v=[a^r]_2\de/\de z^r$ for suitable $a^r\in\ca I_S$ (and we are using local
coordinates adapted to~$S$, as usual). Then from $\tilde X_q(\de_s)=O$ for $s=1,\ldots,m$
it easily follows that $a^r\in\ca I_S^2$ for $r=1,\ldots,m$, and hence $q=O$.

\end{proof}

\begin{remark}
As suggested by \cite{At}, the sequence of the first jets sheaves can be interpreted
as a subsequence of the extension obtained dualizing the sequence
(\ref{eq:extA})\ and tensorizing with~$\ca N_S$:
\[
\def\normalbaselines{\baselineskip20pt \lineskip3pt \lineskiplimit3pt }
\begin{array}{ccccccccc}
O&\longrightarrow&\Omega_S\otimes {\mathcal N}_S&\longrightarrow&
J^1 ({\mathcal N}_S)&\longrightarrow&
{\mathcal N}_S&\longrightarrow& 0\\
&&\Big\|&&&&\Big\downarrow&&\\
O&\longrightarrow& \Omega_S\otimes \ca N_S &\longrightarrow& {\ca A}^*\otimes \ca
N_S&\longrightarrow &\Hom({\ca N}_S,\ca N_S)\otimes
\ca N_S&\longrightarrow& O
\end{array},
\]
where the last vertical map is the injection locally given by
$s\mapsto\id\otimes s$, where $s$ is a local section of ${\ca N}_S$;
it is obtained tensorizing by $ {\ca N}_S$ the inclusion $\ca
O_S\hookrightarrow\Hom({\ca N}_S,\ca N_S)$ that corresponds to the
identity map on  ${\ca N}_S$. In particular, $J^1 ({\mathcal N}_S)$
is a subsheaf of Hom$_{\ca O_S}(\ca A, {\mathcal N}_S)$.
\end{remark}

When $S$ has codimension~1 in~$M$, we have $\Hom({\ca
N}_S,\ca N_S)\cong {\ca O}_S$ and hence the previous remark yields

\begin{corollary}
\label{th:codim1} If $S$ is a codimension~$1$ submanifold of a
complex manifold~$M$, then the sequence $(\ref{eq:extA})$ becomes
$O\to\ca O_S\to{\ca A}\to \ca T_S\to O$, so that $ {\ca A}\cong
J^1({\mathcal N}^*_S) \otimes {\ca N}_S$ and $\ca A$ admits a
nowhere zero holomorphic section.
\end{corollary}

We end this section with a remark that will be useful in Section~\ref{sec:7}:

\begin{proposition}
\label{th:newuno}
Let $S$ be a submanifold of a complex
manifold~$M$, comfortably embedded with respect to a first order lifting~$\rho\colon\ca O_S\to\ca
O_{S(1)}$. Then there exists an $\ca O_S$-morphism
$\tilde\pi\colon\ca T_{M,S(1)}\to
\ca{A}$ such that $\tilde\pi|_{\ca T^S_{M,S(1)}}=\pi$, where $\ca T_{M,S(1)}$ is endowed with the
structure of $\ca O_S$-module given by restriction of scalars via~$\rho$.
\end{proposition}

\begin{proof}
Fix a comfortable atlas $\gt U=\{(U_\alpha,z_\alpha)\}$ adapted to~$\rho$. If
$v\in\ca T_{M,S(1)}$, we can write
\[
v=[x^j_\alpha]_2\frac{\de}{\de z^j_\alpha}
\]
for suitable $[x^j_\alpha]_2\in\ca O_{S(1)}$; then we set
\[
\tilde\pi(v)=\pi\left(\tilde\rho([x^r_\alpha]_2)\frac{\de}{\de z_\alpha^r}+[x^p_\alpha]_2
\frac{\de}{\de z^p_\alpha}\right)=\pi\left(v-\rho([x^r_\alpha]_1)\frac{\de}{\de z_\alpha^r}
\right).
\]
We claim that $\tilde\pi$ is well-defined, that is it does not depend on the particular chart
chosen to express~$v$. Indeed, if we also write $v=[x^k_\beta]_2\frac{\de}{\de z^k_\beta}$
then we have
\[
[x^k_\beta]_2=[x^j_\alpha]_2\left[\frac{\de z^k_\beta}{\de z^j_\alpha}\right]_2.
\]
Applying $\rho\circ\theta_1$ to both sides and recalling that we are working with an
atlas adapted to~$S$ we get
\[
\rho([x^s_\beta]_1)=\rho([x^r_\alpha]_1)\,\rho\left(\left[\frac{\de z^s_\beta}{\de
z^r_\alpha}\right]_1\right).
\]
Then
\[
\rho([x^s_\beta]_1)\frac{\de}{\de z^s_\beta}
=\rho([x^r_\alpha]_1)\,\rho\left(\left[\frac{\de z^s_\beta}{\de
z^r_\alpha}\right]_1\right)\frac{\de}{\de z^s_\beta}
=\rho([x^r_\alpha]_1)\,\left[\frac{\de
z^s_\beta}{\de z^r_\alpha}\right]_2\frac{\de}{\de z^s_\beta}
\]
where we used (\ref{eq:defgrho}) and the fact that $\gt U$ is a comfortable atlas adapted to~$\rho$.
So
\[
\left(v-\rho([x^s_\beta]_1)\frac{\de}{\de z^s_\beta}\right)-\left(v-\rho([x^r_\alpha]_1)
\frac{\de}{\de z^r_\alpha}\right)\in\ca I_S\cdot \ca T^S_{M,S(1)}
\]
because we are using a splitting atlas, and thus $\tilde\pi$ is
well-defined. Finally, it is easy to check that $\tilde\pi$ is
an~$\ca O_S$-morphism extending~$\pi$, and we are done.
\end{proof}

\section{The general index theorem}
\label{sec:6}

In this section we shall prove our general index theorem following the
strategy indicated in the introduction; in the next two sections we shall show how all
the known (and a couple of new ones) index theorems of this kind (with $M$ smooth) for both
holomorphic maps and holomorphic foliations are just particular instances of our general statement.

In the previous sections we have shown how to get a partial holomorphic connection
on the normal bundle from a splitting morphism~$\psi$. The next step is showing how the existence of
a partial holomorphic connection forces the vanishing of some Chern classes. This has been proved,
for instance, by Baum and Bott (\cite{BB}; see also~\cite{CC}); we report here a statement (and a
proof) adapted to our situation.

\begin{theorem}
\label{th:Bott}
Let $S$ be a complex manifold, $F$ a sub-bundle of~$TS$ of rank~$\ell$, and
$E$ a complex vector bundle on~$S$. Assume we have a partial holomorphic connection on~$E$
along~$F$. Then:
\begin{itemize}
\item[\textup{(i)}] every symmetric polynomial in the Chern classes of~$E$ of degree larger than
$\dim S-\ell+\lfloor \ell/2\rfloor$ vanishes.
\item[\textup{(ii)}] Furthermore, if $F$ is involutive and the partial holomorphic connection is
flat then every symmetric polynomial in the Chern classes of~$E$ of degree larger than
$\dim S-\ell$ vanishes.
\end{itemize}
\end{theorem}

\begin{proof}
Write
\begin{equation}
T^{\R}S\otimes\C=F\oplus F_1\oplus T^{(0,1)}S,
\label{eq:decomp}
\end{equation}
where $F_1$ is any $C^\infty$-complement of~$F$
in~$TS=T^{(1,0)}S$. Define a (real) connection~$\nabla$ on~$E$ using the given
partial holomorphic connection on~$F$, any connection on~$F_1$, and $\bar\de$ on~$T^{(0,1)}S$.

Let $\omega$ be the curvature form of~$\nabla$. We claim that
\begin{equation}
\omega(v,\bar{w})=\omega(\bar{u},\bar{w})=O
\label{eq:vc}
\end{equation}
for all $v\in F$ and $\bar{u}$,~$\bar{w}\in T^{(0,1)}S$. It is enough to prove that
they vanish when applied to holomorphic sections of~$E$, since these generate~$\Gamma(E)$ as a
$C^\infty$-module, and the curvature is a tensor. But if $\sigma$ is a holomorphic section of~$E$
we have
\[
\omega(v,\bar{w})(\sigma)=\nabla_v(\nabla_{\bar{w}}\sigma)-
\nabla_{\bar{w}}(\nabla_v\sigma)-\nabla_{[v,\bar{w}]}\sigma=O,
\]
because $\nabla_{\bar{w}}$ kills every holomorphic section, $\nabla_v\sigma$ is
holomorphic because $\nabla$ is holomorphic along~$F$, and $[v,\bar{w}]=O$.
Analogously, since $[\bar{u},\bar{w}]\in T^{(0,1)}S$, one shows that
$\omega(\bar{u},\bar{w})=O$.

Choose local coordinates and local forms~$\eta^1,\ldots,\eta^n$ (where $n=\dim S$) so that
\[
\{\eta^1,\ldots,\eta^\ell,\eta^{\ell+1},\ldots,\eta^n,d\overline{z^1},\ldots,d\overline{z^n}\}
\]
is a local frame for the dual of $T^{\R}S\otimes\C$ respecting~(\ref{eq:decomp}); in particular,
$\{\eta^1|_F,\ldots,\eta^\ell|_F\}$ is a local frame for the dual of~$F$, and
$\{\eta^{\ell+1}|_{F_1},\ldots,\eta^n|_{F_1}\}$ is a local frame for the dual of~$F_1$. Then
(\ref{eq:vc}) implies that in this local frame the curvature matrix is composed by forms which are
linear combinations of
\[
\eta^{p'}\wedge \eta^{q'},\, \eta^{p'}\wedge\eta^{q''},\, \eta^{p''}\wedge\eta^{q''}, \,
d\bar{z^j}\wedge\eta^{q''},
\]
where $1\le p'<q'\le \ell$, $\ell+1\le p''<q''\le n$ and $1\le j\le n$. Since any product of more
than $n-\ell+\lfloor \ell/2\rfloor$ of these forms vanishes, (i) follows.

If $F$ is involutive and the partial holomorphic connection along
$F$ is flat, we moreover have $\omega(v,w)=O$ for all~$v$,~$w\in F$.
This means that we can drop the forms $\eta^{p'}\wedge \eta^{q'}$
from the previous list, and then any product of more than $n-\ell$
of the remaining forms vanishes, giving part (ii).
\end{proof}

\begin{remark}
The previous proof shows not only that Chern classes of suitable degree vanish, but that the
standard differential forms representing them (the one obtained starting from the curvature
matrix of a connection) vanish too.
\end{remark}

We have now all the ingredients needed to apply the general cohomological argument devised by
Lehmann and Suwa. Let us first introduce a couple of definitions to simplify the statements of our
theorems.

\begin{definition}
\label{def:nest}
Let $S$ be a (possibly singular) subvariety of a complex manifold~$M$, and let $S^{\rm reg}\subseteq
S$ be the regular part of~$S$. We shall say that $S$ has an \emph{extendable normal bundle} if there
exists a coherent sheaf of $\ca C^\infty_M$-modules $\ca N$ defined on an open neighborhood of $S$
in $M$ such that $\ca N \otimes_{\ca O_M} \ca
O_{S^{\rm reg}}=\ca N_{S^{\rm reg}}$. We say that $\ca N$ is an \emph{extension} of $\ca
N_{S^{\rm reg}}$.
\end{definition}

\begin{example}
\label{ex:nest}
Any nonsingular submanifold has an extendable normal bundle:
an extension of~$\ca N_S$ is given by the pull-back (under the retraction)
to a tubular neighbourhood. If $S$ is singular but has codimension one in~$M$
then the line bundle $\ca O([S])$ associated to the divisor $[S]$
provides an extension of~$\ca N_{S^{\rm reg}}$, and hence $S$ has an extendable normal bundle. More
generally, if $S$ is a \emph{locally complete intersection defined by
a section,} or a \emph{strongly locally complete intersection,}
then it has an extendable normal bundle (see \cite{LS1}
and~\cite{LS2}).
\end{example}

\begin{remark}
\label{rem:nest}
The extension of the normal bundle might be, in general, not unique. However, in all cases
described in the previous example there is a natural extension to consider.
\end{remark}

The next definition will considerably shorten several statements.

\begin{definition}
Let $S$ be a compact, complex, reduced, irreducible, possibly singular, subvariety of
dimension~$d$ of an $n$-dimensional complex manifold~$M$. Assume that $S$ has extendable normal
bundle. Let $\Sigma$ be an analytic subset of~$S$, containing the singular
part~$S^{\rm{sing}}$ of $S$, so that $S^o=S\setminus\Sigma\subseteq S^{\rm reg}$, and let
furthermore
$\gt F$ denote another analytic object involved in the problem (for instance, in our applications
$\gt F$ will be either a holomorphic foliation or a holomorphic self-map, and $\Sigma$ the union of
the singular set of~$S$ with the singular set of~$\gt F$). We shall say that $S$
\emph{has the Lehmann-Suwa index property of level~$\ell\ge 1$ on~$\Sigma$ with respect to~$\gt F$}
if given an extension $\ca N$ of~$\ca N_{S^{\rm reg}}$ we can associate to every homogeneous
symmetric polynomial~$\phe$ of degree~$k>d-\ell$ and every connected component~$\Sigma_\lambda$
of~$\Sigma$ a homology class
\[
\hbox{\rm Res}_{\phe}(\gt F,\ca N;\Sigma_\lambda)\in
H_{2(d-k)}(\Sigma_\lambda;\C),
\]
depending only on $\ca N$ and on the local behavior of~$\gt F$
near~$\Sigma_\lambda$, so that
\[
\sum (i_\lambda)_*\hbox{\rm Res}_{\phe}(\gt F,\ca N;\Sigma_\lambda)=[S]\frown\phe(\ca N)
\quad\hbox{in $H_{2(d-k)}(S;\C)$,}
\]
where the sum ranges over all the connected components of~$\Sigma$, the map $i_\lambda\colon
\Sigma_\lambda\hookrightarrow S$ is the inclusion, and $\phe(\ca N)$ denotes the class obtained
evaluating~$\phe$ in the Chern classes of~$\ca N$.
\end{definition}

\begin{remark}
When $k=d$ then $[S]\frown\phe(\ca N)=\int_S\phe(\ca N)\in\C$.
\end{remark}

And now, our general index theorem:

\begin{theorem}
\label{th:generalindex}
Let $S$ be a compact, complex, reduced, irreducible, possibly singular, subvariety of
dimension~$d$ of an $n$-dimensional complex manifold~$M$, and assume that $S$ has extendable normal
bundle. Let $\Sigma$ be an analytic subset of $S$ containing~$S^{\rm sing}$ such that there exist a
sub-bundle $F$ of rank~$\ell$ of~$TS^o$ (where $S^o=S\setminus\Sigma\subseteq S^{\rm reg}$) and an
$\ca O_{S^o}$-morphism $\psi\colon\ca F\to\ca A$ with $\theta_1\circ\psi=\id$, where $\ca A$ is the
Atiyah sheaf of~$S^o$.  Then:
\begin{itemize}
\item[\textup{(i)}] $S$ has the Lehmann-Suwa index property of level~$\ell-\lfloor\ell/2\rfloor$
on~$\Sigma$ with respect to~$\psi$.
\item[\textup{(ii)}] If furthermore $F$ is involutive and $\psi$ is a Lie algebroid morphism, then
$S$ has the Lehmann-Suwa index property of level~$\ell$ on~$\Sigma$ with respect
to~$\psi$.
\end{itemize}
\end{theorem}

\begin{proof}
Theorem~\ref{th:sommario} yields a partial holomorphic connection
on~$N_{S^o}$ along~$F$, which is flat in case (ii).
Theorem~\ref{th:Bott} then implies that every symmetric polynomial
in the Chern classes of~$N_{S^o}$ of degree larger than
$d-\ell+\lfloor \ell/2\rfloor$ (or, in case (ii), larger
than~$d-\ell$) vanishes. The theorem then follows from the general
Lehmann-Suwa theory (see, e.g., Chapter V\negthinspace I in
\cite{S2}, or \cite{LS2}).
\end{proof}

\begin{remark}
For the sake of completeness, let us summarize here the gist of Lehmann-Suwa's argument, warning
the reader that the complete proof is a bit technical and requires \v Cech-de Rham
cohomology (see, e.g., \cite{S2}). Let
\[
H^*(S,S^o;\C)\longrightarrow H^*(S;\C)\longrightarrow H^*(S^o;\C)
\]
be the long exact cohomology sequence of the pair~$(S,S^o)$. The vanishing Theorem~\ref{th:Bott}
says that the cohomology class~$\phe(\ca N)$ vanishes when restricted to~$S^o$; hence it must be
the image of some cohomology class~$\eta\in H^*(S,S^o;\C)$ which is, by definition of cohomology of
a pair, concentrated in an arbitrary neighbourhood of~$S\setminus S^o=\Sigma$. Such a class is not
unique in general, and it should be chosen in a suitable way depending on the partial holomorphic
connection given by~$\psi$. Now, since
$S$ is compact, the Poincar\'e homomorphism (consisting exactly in taking the cap product
with~$[S]$) gives a natural map from $H^*(S;\C)$ to~$H_{2d-*}(S;\C)$. On the other
hand, since
$\Sigma$ is an analytic subset of~$S$, the Alexander homomorphism~$A$ gives a natural map
from~$H^*(S,S^o;\C)$ to~$H_{2d-*}(\Sigma;\C)$. Furthermore, if we denote by $i\colon\Sigma\to S$ the
inclusion, we have the equality $i_*A(\eta)=[S]\frown\phe(\ca N)$. Now, if $\Sigma=\bigcup_\lambda
\Sigma_\lambda$ is the decomposition in connected components of~$\Sigma$, we have
$H_{2d-*}(\Sigma;\C)=\bigoplus_\lambda H_{2d-*}(\Sigma_\lambda;\C)$; therefore if we denote by
$\hbox{\rm Res}_\phe(\psi,\ca N;\Sigma_\lambda)$ the component of~$A(\eta)$ belonging
to~$H_{2d-*}(\Sigma_\lambda;\C)$, we obtain
\[
\sum_\lambda (i_\lambda)_* \hbox{\rm Res}_\phe(\psi,\ca N;\Sigma_\lambda)=[S]\frown\phe(\ca N),
\]
that is the index theorem.
\end{remark}

\begin{remark}
\label{rem:Filippo1}
Let us now describe how to compute $\hbox{\rm Res}_\phe(\psi;\Sigma_\lambda)$ in a simple (but
useful) case. Assume that: $S$ is a locally
complete intersection defined by a section; the connected component~$\Sigma_\lambda$ reduces to an
isolated point~$p\in S$; the sub-bundle~$F$ has rank~$\ell=1$; and there exists a local
vector field $v\in(\ca T_S)_p\subset(\ca T_{M,S})_p$ vanishing at~$p$ and generating~$F$ in a
pointed neighbourhood of~$p$. Let $l^1,\ldots,l^m$ be a local system of defininig functions for~$S$
near~$p$, so that $\{[l^1]_2,\ldots,[l^m]_2\}$ is a local frame for $\ca I_S/\ca I_S^2=\ca N_S^*$,
and denote by $\{\zeta_1,\ldots,\zeta_m\}$ the corresponding dual local frame of~$\ca N_S$. If
$\tilde X^\psi$ is the partial holomorphic connection induced by~$\psi$, then writing
\[
\tilde X^\psi_v(\zeta_r)=c^s_r\zeta_s
\]
we get an $m\times m$ matrix $C=(c^r_s)$ of holomorphic functions defined in a pointed
neighbourhood of~$p$. Finally (see~\cite{LS1}), it is possible to choose a local chart $(U,z)$
at~$p$ so that if we write $v=[a^j]_1\de/\de z^j$ then
\[
\{l^1=\cdots=l^m=a^{m+1}=\cdots=a^n=0\}=\{p\}
\]
(if $p$ is a regular point of~$S$ it suffices to take any chart adapted to~$S$). Take now a
homogeneous symmetric polynomial~$\phe$ of degree $\dim S$; then $\hbox{\rm Res}_\phe(\psi;\{p\})$
is given by the Grothendieck residue
\begin{equation}
\hbox{\rm Res}_\phe(\psi,\ca N;\{p\})=\frac{1}{(2\pi
i)^{n-m}}\int_\Gamma\frac{\phe(C)}{a^{m+1}\cdots a^n}\,dz^{m+1}\wedge
\cdots\wedge dz^n,
\label{eq:Grot}
\end{equation}
where $\Gamma=\{q\in S\mid |a^{m+1}(q)|=\cdots=|a^n(q)|=\varepsilon\}$ for $0<\varepsilon<<1$,
oriented so that $d\arg a^{m+1}\wedge\cdots\wedge d\arg a^n$ is positive, $\ca N$ is the natural
extension of~$\ca N_{S^{\rm reg}}$ mentioned in Remark~\ref{rem:nest}, and
$\phe(C)$ denotes~$\phe$ evaluated on the eigenvalues of the matrix~$C$. This formula can be
obtained by observing that if $\tilde v\in\ca T^S_{M,S(1)}$ is such that $\pi(\tilde v)=\psi(v)$
outside~$p$, then the local partial holomorphic connection on~$\ca N_{S^{\rm reg}}$ induced
by~$\tilde v$ coincides with~$\tilde X^\psi_v$, and the residue $\hbox{\rm Res}_\phe(\psi,\ca
N;\{p\})$ coincides with the residue associated to~$\tilde v$ and obtained in~\cite{LS1}.
By the way, an explicit algorithm for computing the Grothendieck residue~(\ref{eq:Grot}) when
$p$ is a regular point of~$S$ is described in~\cite{BB}, p. 280.
\end{remark}

We end this section by describing the general strategy we are going to use to build the
morphism~$\psi\colon\ca F\to\ca A$. Such a morphism exists if and only if the sequence
\[
O\longrightarrow\Hom(\ca N_{S^o},\ca N_{S^o})\longrightarrow\theta_1^{-1}(\ca F)\stackrel{\theta_1}
{\longrightarrow}\ca F\longrightarrow O
\]
splits, that is if and only if the associated cohomology class in $H^1\bigl(S^o,\ca F^*\otimes
\Hom(\ca N_{S^o},\ca N_{S^o})\bigr)$ vanishes. The latter class is represented by a cocycle of the
form~$\{\psi_\beta-\psi_\alpha\}$, where the $\psi_\alpha$ are local splitting morphisms. Therefore
the morphism~$\psi$ exists if and only if we can find local morphisms~$x_\alpha$ from~$\ca F$
to~$\Hom(\ca N_{S^o},\ca N_{S^o})$ such that $\psi_\beta-\psi_\alpha=x_\beta-x_\alpha$.

Our strategy then will be to use the additional data involved (foliation or self-map) to build
local splitting morphisms; in this way we shall be able to express the cohomological problem in
terms of the geometry of the additional data, and then to give sufficient
conditions for the problem to be solvable.

Notice in particular that if $S^o$ is Stein then this cohomological problem is \emph{always}
solvable, and thus we have

\begin{corollary}
\label{th:Steinindex}
Let $S$ be a compact, complex, reduced, irreducible, possibly singular, subvariety of
dimension~$d$ of an $n$-dimensional complex manifold~$M$, and assume that $S$ has extendable normal
bundle. Let $\Sigma$ be an analytic subset containing~$S^{\rm sing}$ such that
$S^o=S\setminus\Sigma$ is Stein.  Then $S$ has the Lehmann-Suwa index property of level~$d-\lfloor
d/2\rfloor$ on~$\Sigma$ with respect to anything providing local splitting morphisms for the
sequence~$(\ref{eq:extA})$ over~$S^o$.
\end{corollary}

\section{Holomorphic foliations}
\label{sec:7}

The aim of this section is to show how to use a holomorphic foliation on the ambient
manifold to implement the strategy just discussed. We recall that a (possibly singular)
\emph{holomorphic foliation~$\ca F$} of \emph{dimension}~$\ell$ on a complex $n$-dimensional
manifold~$M$ is, by definition, a coherent involutive subsheaf of~$\ca T_M$ which is a locally
free $\ca O_M$-module of rank~$\ell$ outside its singular locus~$\Sing(\ca F)$, defined as the set
of points
$x\in M$ such that the quotient $\ca T_M/\ca F$ is not locally free of dimension~$n-\ell$ at~$x$.
The foliation is called
\emph{non-singular} if the singular locus is empty. We refer to \cite{BB} and \cite{S2},
chapter~V\negthinspace I, for more details on holomorphic foliations.

\begin{definition} Let $S$ be a complex (not necessarily closed) $m$-codimensional submanifold of
an $n$-dimensional complex manifold~$M$, and let $\ca F$ be a (possibly singular) holomorphic
foliation~$\ca F$  on~$M$, of dimension $\ell\le n-m=\dim S$. We shall denote
by~$\ca F_{S(1)}$ the $\ca O_{S(1)}$-submodule~$\ca F\otimes_{\ca O_M}\ca O_{S(1)}$ of~$\ca
T_{M,S(1)}$, and by~$\ca F_S$ the $\ca O_S$-submodule $\ca F\otimes_{\ca O_M}\ca O_S$ of~$\ca
T_{M,S}$. If $\ca F_S\subseteq\ca T_S\subset
\ca T_{M,S}$, then $\ca F$ is \emph{tangent} to $S$; otherwise, the foliation $\ca F$ is
\emph{tranverse} to
$S$.
\end{definition}

\begin{remark}
We shall always assume that $S$ is not contained in the singular locus of~$\ca F$.
\end{remark}

In the tangential case, we clearly have $\ca F_{S(1)}\subseteq\ca T^S_{M,S(1)}$.
Furthermore, $\ca F_S$ is a (possibly singular) holomorphic foliation of~$S$ of dimension~$\ell$.
The singular locus of~$\ca F_S$ (which is the intersection of~$\Sing(\ca F)$ with~$S$) is an
analytic subset of~$S$; therefore since our aim is to build a splitting morphism $\psi$
outside the singularities, we shall assume that

\begin{case}
$\ca F_S$ is a non-singular holomorphic foliation of~$S$ of dimension~$\ell\le\dim S$ (and thus, in
particular, it is the sheaf of germs of holomorphic sections of an involutive sub-bundle~$F$
of~$TS$ of rank~$\ell$). To be consistent with the non-tangential case, we shall also set
$\ca F^\sigma=\ca F_S$ and~$\sigma^*=\id_{\ca F_S}$.
\end{case}

If $\ca F$ is not tangent to~$S$ then $\ca F_S$ is not a subsheaf of~$\ca T_S$, but only of~$\ca
T_{M,S}$. To get a subsheaf of~$\ca T_S$, we must project~$\ca F_S$ into it.

\begin{definition} Let $S$ be a splitting submanifold of a complex manifold~$M$. Given a first order
lifting~$\rho\colon\ca O_S\to\ca O_M/\ca I_S^2$, let $\sigma^*\colon\ca T_{M,S}\to\ca T_S$
be the left splitting morphism associated to~$\rho$ by Proposition~\ref{th:ucinque}. If $\ca F$ is a
holomorphic foliation on~$M$ of dimension~$\ell\le\dim S$, we shall denote by~$\ca F^\sigma$ the coherent
sheaf of
$\ca O_S$-modules given by
\[
\ca F^\sigma=\sigma^*(\ca F_S)\subseteq\ca T_S.
\]
We shall say that $\rho$ is \emph{$\ca F$-faithful outside an analytic subset~$\Sigma\subset S$}
if $\ca F^\sigma$ is a non-singular holomorphic foliation of dimension~$\ell$ on~$S\setminus\Sigma$.
If $\Sigma=\void$ we shall simply say that $\rho$ is \emph{$\ca F$-faithful.}
\end{definition}

It might happen that a first order lifting is not $\ca F$-faithful
while another one is. Furthermore, $\ca F^\sigma$ might be as well as
not be involutive depending on the choice of~$\rho$.

\begin{example}
Let $M=\C^4$, take $S=\{z^1=0\}$ and let $\ca F$ be the
non-singular foliation generated over~$\ca O_M$ by the global vector
fields $(z^2-z^1)\frac{\de}{\de z^3}+\frac{\de}{\de z^4}$ and
$\frac{\de}{\de z^1}+\frac{\de}{\de z^2}$, so that $\ca F_S$ is
generated over~$\ca O_S$ by~$z^2\frac{\de}{\de z^3}+\frac{\de}{\de z^4}$
and $\frac{\de}{\de z^1}+\frac{\de}{\de z^2}$. The submanifold $S$
clearly splits in~$M$, and a natural choice of first order lifting
is
\[
\rho([f]_1)=[f]_2-\left[\frac{\de f}{\de z^1}z^1\right]_2.
\]
The corresponding left splitting morphism $\sigma^*$ is the identity on~$\ca T_S$ and kills
$\frac{\de}{\de z^1}$; therefore $\ca F^\sigma$ is generated over~$\ca
O_S$ by $z^2\frac{\de}{\de z^3}+\frac{\de}{\de z^4}$ and $\frac{\de}{\de
z^2}$, and thus $\ca F^\sigma$ is not involutive.

If we choose as first order lifting the less standard $\rho_1$
given by
\[
\rho_1([f]_1)=[f]_2-\left[\left(\frac{\de f}{\de z^1}+\frac{\de f}{\de z^2}\right)z^1
\right]_2,
\]
then the corresponding left splitting morphism~$\sigma_1^*$ sends $\frac{\de}{\de z^1}$
in~$-\frac{\de}{\de z^2}$, and $\ca F^{\sigma_1}$ turns out to be generated by~$z^2\frac{\de}{\de
z^3}+\frac{\de}{\de z^4}$ only, and so it is involutive, but of the wrong dimension.
Finally, if we take as first order lifting
\[
\rho_2([f]_1)=[f]_2-\left[\left(\frac{\de f}{\de z^1}+\frac{\de f}{\de z^2}-\frac{\de f}{\de
z^3}\right)z^1\right]_2
\]
then $\sigma_2^*$ sends $\frac{\de}{\de z^1}$
in~$\frac{\de}{\de z^3}-\frac{\de}{\de z^2}$, so that $\ca F^{\sigma_2}$ is generated
by~$z^2\frac{\de}{\de z^3}+\frac{\de}{\de z^4}$ and $\frac{\de}{\de z^3}$, and thus
$\rho_2$ is $\ca F$-faithful.
\end{example}

If $\ca F^\sigma$ has dimension equal to~1 or to the dimension of~$S$, then it is
automatically involutive. In this case it is easy to have faithfulness:

\begin{lemma}
\label{th:dfaith}
Let $S$ be a splitting submanifold of a complex manifold~$M$, and let $\ca F$ be a
holomorphic foliation on~$M$ of dimension equal to~$1$ or to the dimension of~$S$. If there exists $x_0
\in S\setminus\Sing(\ca F)$ such that $\ca F$ is tangent
to $S$ at $x_0$, i.e., $(\ca F_S)_{x_0} \subseteq \ca T_{S,x_0}$,
then any first order lifting is $\ca F$-faithful outside a suitable analytic subset of~$S$.
\end{lemma}

\begin{proof}
Let $\sigma^*\colon\ca T_{M,S}\to\ca T_S$ be the left-splitting
morphism associated to a first order lifting~$\rho$. By assumption,
$\ker\sigma^*_{x_0}\cap\ca (F_S)_{x_0}=(O)$; therefore
$\ker\sigma^*_x\cap\ca (F_S)_x=(O)$ for all $x\in S$ outside an
analytic subset~$\Sigma_0$ of~$S$. Furthermore, $\ca F^\sigma$ has
dimension equal to~1 or to~$\dim S$; therefore it is involutive, and
hence $\rho$ is $\ca F$-faithful outside~$\Sigma_0\cup\Sing(\ca
F^\sigma)$.
\end{proof}

Notice that there are topological obstructions for a foliation to
be everywhere non-tangential to $S$. For instance, in \cite{Bru}, \cite{Hon}
it is proved that if $S$ is a curve in a surface $M$, the number of
points of tangency between $S$ and an one-dimensional holomorphic reduced foliation~$\ca F$ of~$M$,
counted with multiplicity, is $S\cdot S-S\cdot\ca F$. Therefore if $S\cdot S\ne S\cdot\ca F$ then
every first order lifting is $\ca F$-faithful outside a suitable analytic subset.

Another result of this kind shows that for one-dimensional foliations most first order liftings
are faithful:

\begin{lemma}
\label{th:dsplitfaith}
Let $S$ be a non-singular hypersurface splitting in a
complex manifold~$M$, and let $\ca F$ be a one dimensional
holomorphic foliation on~$M$. Assume that $S$ is not contained in
$\Sing(\ca F)$. Then there is at most one first order lifting
$\rho$ which is not $\ca F$-faithful outside a suitable analytic subset of~$S$.
\end{lemma}

\begin{proof} Suppose $\rho$ is a first order lifting of $S$ which
is \emph{not} $\ca F$-faithful; since $\ca F$ is one-dimensional,
this means that $(\ca F_S)_x\subseteq\ker\sigma^*_x$ for all $x\in
S\setminus\Sing(\ca F)$, where $\sigma^*$ is the left splitting
morphism associated to~$\rho$. By Lemma~\ref{th:uquattro}.(iii) any
other left splitting morphism is of the form
$\sigma_1^*=\sigma^*+\phe \circ p_2$ with~$\phe \in
H^0\bigl(S,\Hom(\ca N_S,\ca T_S)\bigr)$; in particular,
$\sigma_1^*(v)=\phe\bigl(p_2(v)\bigr)$ for all~$v\in\ca F_S$. Now,
since $\sigma^*$ is a left splitting morphism, we have
$\ker\sigma^*_x\cap\ker (p_2)_x=(O)$ for all~$x\in S$; therefore
$p_2|_{\ca F_S}$ is injective. Furthermore, since $\ca N_S$ has rank
one, $\phe_x$ is either injective or identically zero; hence
$(\sigma_1^*)_x$ restricted to~$(\ca F_S)_x$ is either injective or
identically zero. Now, if $\phe\ne O$ then $\phe_x\ne O$ for $x$
outside an analytic subset~$\Sigma_0$ of~$S$; therefore it follows
that if $\phe\ne O$ then the first order lifting associated to
$\sigma_1^*$ is $\ca F$-faithful outside $\Sigma_0\cup\Sing(\ca
F^{\sigma_1})$.
\end{proof}

\begin{corollary}
\label{th:dfaithalways}
Let $S$ be a non-singular hypersurface splitting in a complex manifold~$M$, and let
$\ca F$ be a one dimensional holomorphic foliation on~$M$. Assume that $S$ is not
contained in~$\Sing(\ca F)$. If $H^0(S,\ca T_S \otimes \ca
N_S^*)\neq (O)$ then there exists at least one first order lifting $\ca F$-faithful
outside a suitable anlytic subset of~$S$.
\end{corollary}

\begin{proof}
If $H^0(S,\ca T_S \otimes \ca N_S^*)\neq(O)$ then by
Lemma~\ref{th:uquattro}.(iii) there exist at least two different
splitting morphisms. Then the assertion follows from the previous
lemma.
\end{proof}

Coming back to our main concern, in the non-tangential case we shall momentarily make the following
assumption:

\begin{case} There exists an $\ca F$-faithful first order lifting $\rho$, with associated
left-splitting morphism~$\sigma^*$; in particular, $\ca F^\sigma$ is a non-singular holomorphic
foliation of~$S$ of dimension~$\ell\le\dim S$, and $\sigma^*|_{\ca F_S}\colon\ca F_S\to\ca
F^\sigma$ is an isomorphism of~$\ca O_S$-modules.
\end{case}

If $\ca G\subseteq\ca T_S$ is a non-singular holomorphic foliation
of~$S$ of dimension~$\ell\le\dim S$, Frobenius' theorem implies that we can always find an
atlas~$\gt U=\{(U_\alpha,z_\alpha)\}$ adapted to~$S$ such that the $\{\de/\de
z^{m+1}_\alpha,\ldots,\de/\de z^{m+\ell}_\alpha\}$ are local frames
for~$\ca G$. Furthermore, it is easy to check that if $S$ is split
(2-split, comfortably embedded) in~$M$ we can also assume that $\gt U$
is a splitting (2-splitting, comfortable) atlas.

\begin{definition}
Assume we are either in Case~1 or in Case~2. An atlas~$\gt
U=\{(U_\alpha,z_\alpha)\}$ adapted to~$S$ such that the $\{\de/\de
z^{m+1}_\alpha,\ldots,\de/\de z^{m+\ell}_\alpha\}$ are local frames
for~$\ca F^\sigma$ shall be said \emph{adapted to~$S$ and~$\ca F$.} We
explicitly notice that if $\gt U$ is adapted to $S$ and~$\ca F$ then
\[
\frac{\de z_\beta^{q''}}{\de z_\alpha^{p'}}\in\ca I_S
\]
for all $p'=m+1,\ldots,m+\ell$, $q''=m+\ell+1,\ldots,n$ and indices $\alpha$, $\beta$ such that
$U_\alpha\cap U_\beta\cap S\ne\void$.
\end{definition}

\begin{remark}
From now on, indices like $p'$, $q'$, $\tilde p'$ and $\tilde q'$ will run from $m+1$
to~$m+\ell$, while indices like~$p''$ and~$q''$ will run from~$m+\ell+1$ to~$n$.
\end{remark}

Using adapted atlas we can find special local frames for the foliation:

\begin{lemma}
\label{th:newinter}
Assume we are in Case~1 or in Case~2, and let
$\{(U_\alpha,z_\alpha)\}$ be an atlas adapted to~$S$ and~$\ca F$ (and to $\rho$ too in Case~2).
\begin{itemize}
\item[\textup{(i)}] For each index
$\alpha$ there exists a unique $\ell$-uple~$(v_{\alpha,m+1},\ldots, v_{\alpha,m+\ell})$
of elements of~$\ca F$ of the form
\begin{equation}
v_{\alpha,p'}=\frac{\de}{\de z_\alpha^{p'}}+(a_\alpha)^r_{p'}\frac{\de}{\de
z_\alpha^r}+(a_\alpha)^{p''}_{p'}\frac{\de}{\de
z_\alpha^{p''}}
\label{eq:newform}
\end{equation}
with $(a_\alpha)^{p''}_{p'}\in\ca I_S$
for $p'=m+1,\ldots,m+\ell$ and~$p''=m+\ell+1,\ldots,n$,
(and, in Case~1,
$(a_\alpha)^r_{p'}\in\ca I_S$ for~$r=1,\ldots,m$), so that
$\sigma^*(v_{\alpha,p'}\otimes[1]_1)=\de/\de z_\alpha^{p'}$.
\item[\textup{(ii)}] The set $\{v_{\alpha,m+1},\ldots, v_{\alpha,m+\ell}\}$ is a local frame for the
sheaf $\ca F$ in a neighbourhood of $S\cap U_\alpha$. Furthermore, writing
\[
v_{\beta,q'}=(c_{\beta\alpha})^{p'}_{q'} v_{\alpha,p'}
\]
the $(c_{\beta\alpha})^{p'}_{q'}\in\ca O_M$ define a cocycle $(c_{\beta\alpha})$
representing the vector bundle associated to $\ca F$ in a neighbourhood of $S$ and satisfy
the following relations:
\begin{equation}
\begin{cases}\displaystyle (c_{\beta\alpha})^{p'}_{q'}= \frac{\de
z_\alpha^{p'}}{\de z_\beta^{q'}}+ (a_\beta)^s_{q'}\frac{\de
z_\alpha^{p'}}{\de z_\beta^s} + (a_\beta)^{q''}_{q'}\frac{\de
z_\alpha^{p'}}{\de z_\beta^{q''}} ,\cr \displaystyle
(c_{\beta\alpha})^{p'}_{q'} (a_\alpha)^r_{p'}= \frac{\de
z_\alpha^r}{\de z_\beta^{q'}}+ (a_\beta)^s_{q'}\frac{\de
z_\alpha^r}{\de z_\beta^s} +(a_\beta)^{q''}_{q'}\frac{\de
z_\alpha^r}{\de z_\beta^{q''}},\cr \displaystyle
(c_{\beta\alpha})^{p'}_{q'} (a_\alpha)^{p''}_{p'}= \frac{\de
z_\alpha^{p''}}{\de z_\beta^{q'}} + (a_\beta)^s_{q'}\frac{\de
z_\alpha^{p''}}{\de z_\beta^s} +(a_\beta)^{q''}_{q'}\frac{\de
z_\alpha^{p''}}{\de z_\beta^{q''}}. \end{cases} \label{eq:calfabeta}
\end{equation}
\end{itemize}
\end{lemma}

\begin{proof}
(i) Since $\{(U_\alpha,z_\alpha)\}$ is adapted to~$S$ and~$\ca F$, the $\de/\de
z_\alpha^{p'}$'s form local frames for~$\ca F^\sigma$. Hence we can find~$\tilde
v_{\alpha,p'}\in\ca F$ such that
$\sigma^*(\tilde v_{\alpha,p'}\otimes[1]_1)=\de/\de z_\alpha^{p'}$. Write
\[
\tilde v_{\alpha,p'}=(b_\alpha)^{q'}_{p'}\frac{\de}{\de
z_\alpha^{q'}}+(b_\alpha)^r_{p'}\frac{\de}{\de
z_\alpha^r}+(b_\alpha)^{q''}_{p'}\frac{\de}{\de z_\alpha^{q''}}
\]
for suitable $(b_\alpha)^j_{p'}\in\ca O_M$; we must have
$[(b_\alpha)^q_{p'}]_1=\delta^q_{p'}$, and $[(b_\alpha)^r_{p'}]_1=0$ in Case~1. In
particular,
$([(b_\alpha)^{q'}_{p'}]_1)$ is the identity matrix; hence
$\bigl((b_\alpha)^{q'}_{p'}\bigr)$ is invertible as matrix of germs. Multiplying then the
$\ell$-uple~$(\tilde v_{\alpha,m+1},\ldots,\tilde v_{\alpha,m+\ell})$ by the inverse of this
matrix we get an $\ell$-uple $(v_{\alpha,m+1},\ldots, v_{\alpha,m+\ell})$ of elements of~$\ca F$ of
the desired form. Furthermore, since ${\rm rk}_{\ca O_M}\ca F=\ell$, the $v_{\alpha,p'}$'s form a
local frame for $\ca F$; an exercise in linear algebra then shows that they are uniquely
determined.
\smallskip

(ii) The elements $v_{\alpha,m+1},\ldots, v_{\alpha,m+\ell}$ form a
local frame for the sheaf $\ca F$  in a neighbourhood of $S\cap
U_\alpha$ since their restriction to $S$ form a local frame for $\ca
F_S$ on $U_\alpha\cap S$. The relations (\ref{eq:calfabeta})  then
follows directly from (\ref{eq:newform}).
\end{proof}

Restricting to $S$ the first equation in (\ref{eq:calfabeta}) we get:

\begin{corollary}
\label{th:isoF}
Assume we are in Case~1 or in Case~2, let
$\{(U_\alpha,z_\alpha)\}$ be an atlas adapted to~$S$ and~$\ca F$ (and to $\rho$ too in Case~2),
and let $v_{\alpha,m+1},\ldots,v_{\alpha,m+\ell}$ be given by the previous lemma.
Then the vector bundles associated to ${\ca F}_S$ and ${\ca
F^\sigma}$ are represented by the same cocycle
\begin{equation}
[(c_{\beta\alpha})^{p'}_{q'}]_1=
\left[\frac{\de z_\alpha^{p'}}{\de z_\beta^{q'}}\right]_1.
\label{eq:yae}
\end{equation}
in the frames
$\{v_{\alpha,m+1}\otimes[1]_1,\ldots, v_{\alpha,m+\ell}\otimes[1]_1\}$ and $\{\de/\de
z_\alpha^{m+1},\ldots,\de/\de z_\alpha^{m+\ell}\}$, respectively. In particular,
the isomorphism $\sigma^*|_{\ca F_S}\colon {\ca F}_S\to{\ca F^\sigma}$
is represented with respect to these frames by the identity matrix.
\end{corollary}

The restriction to $S$ of the second equation in (\ref{eq:calfabeta}) gives
\[
[(c_{\beta\alpha})^{p'}_{q'}]_1 [(a_\alpha)^r_{p'}]_1
=[(a_\beta)^s_{q'}]_1\left[\frac{\de z_\alpha^r}{\de
z_\beta^s}\right]_1,
\]
and thus we get a global section~$T\in H^0(S,(\ca F^\sigma)^*\otimes\ca N_S)$ by setting
\[
T|_{U_\alpha}=[(a_\alpha)^r_{p'}]_1\,\omega_\alpha^{p'}\otimes\de_{r,\alpha},
\]
where $\omega_\alpha^{p'}$ is the local section of~$(\ca F^\sigma)^*$ induced
by~$dz_\alpha^{p'}$. It is easy to check that the corresponding morphism
$T\colon\ca F^\sigma\to\ca N_S$ is given by the composition
\[
 T\colon \ca F^\sigma \stackrel{(\sigma^*|_{\ca F_S})^{-1}}{\longrightarrow} \ca
F_S\longhookrightarrow
\ca T_{M,S}\stackrel{p_2}{\longrightarrow}
 \ca N_S.
\]

\begin{remark}
\label{rem:transfer}
The morphism $T$ is non-zero if and only if the foliation ${\ca F}$ is
transversal to $S$.
\end{remark}

Now we are ready to characterize the existence of
morphisms $\psi\colon \ca F^\sigma\to {\ca A}$ such that $\theta_1\circ \psi = \id$:

\begin{proposition}
\label{th:generalpsi}
Assume we are in Case~1, or in Case~2 with $S$ comfortably embedded in~$M$.
Given a (comfortable in Case~2) atlas $\gt U=\{(U_\alpha,z_\alpha)\}$ adapted to~$S$ and~$\ca F$
(and to $\rho$ too in Case~2), let $\{f_{\beta\alpha}\}$ be the cocycle defined by
\begin{eqnarray*}
f_{\beta\alpha}&=&
\left[(c_{\alpha\beta})^{q'}_{p'}\right]_1
\tilde\rho\left(\left[(c_{\beta\alpha})^{\tilde p'}_{q'}\right]_2\right)
\otimes\frac{\de}{\de z_\alpha^{\tilde p'}}\otimes\omega_\alpha^{p'}\\
&=&
\left.\frac{\de z_\beta^{q'}}{\de z_\alpha^{p'}}
 \frac{\de(c_{\beta
\alpha})^{\tilde p'}_{q'}}{\de z^t_\alpha}\right|_S
\left[z^t_\alpha\right]_2
\otimes\frac{\de}{\de z_\alpha^{\tilde p'}}\otimes\omega_\alpha^{p'},
\end{eqnarray*}
and denote by $\gt f \in H^1\bigl(S,\ca N_S^*\otimes\ca F^\sigma\otimes(\ca F^\sigma)^*\bigr)$ the
corresponding cohomology class. Then there exists a morphism
$\psi\colon{\ca F^\sigma}\to {\ca A}$ such that $\theta_1\circ \psi=\id$ if and
only if
\[
T_*(\gt f)=O
\]
in $H^1\bigl(S,\ca N_S^*\otimes \ca N_S\otimes(\ca F^\sigma)^*\bigr)$, where $T_*$ is
the map induced in cohomology by the morphism $\id\otimes T\otimes\id$.
\end{proposition}

\begin{proof}
It is easy to check that $\gt f$ is a well-defined cohomology class
indepedent of the particular atlas $\gt U$ chosen. Thus to prove the
assertion it suffices to find local splittings $\psi_\alpha\colon
\ca F^\sigma|_{U_\alpha\cap S}\to\ca A|_{U_\alpha\cap S}$ so that
$\{\psi_\beta-\psi_\alpha\}$ represents the cohomology class~$T_*(\gt f)$.

A local frame for~$\ca F^\sigma$ is~$\{\de/\de z^{m+1}_\alpha,\ldots,\de/\de z^{m+\ell}_\alpha\}$.
We then define $\psi_\alpha$ by setting
\begin{equation}
\psi_\alpha\left(\frac{\de}{\de z^{p'}_\alpha}\right)=\tilde\pi(v_{\alpha,p'}\otimes[1]_2)=
\pi\left(\frac{\de}{\de z_\alpha^{p'}}+
\tilde{\rho}([(a_\alpha)^r_{p'}]_2)
 \frac{\de}{\de z_\alpha^{r}}\right)
\label{eq:psialpha}
\end{equation}
and then extending by $\ca O_S$-linearity, where $\pi\colon\ca T^S_{M,S(1)}\to\ca A$ is
the canonical projection, and $\tilde\pi\colon\ca T_{M,S(1)}\to\ca A$ is the morphism
introduced in Proposition~\ref{th:newuno}. Notice that, in Case~1,
$v_{\alpha,p'}\otimes[1]_2\in\ca T^S_{M,S(1)}$, and so $\psi_\alpha$ is defined without
assuming anything on the embedding of~$S$ into~$M$.

Now, we have:
\begin{eqnarray*}
\psi_\beta\left(\frac{\de}{\de z^{p'}_\alpha}\right)
&=&\psi_\beta\left(\left[\frac{\de z_\beta^{q'}}{\de z_\alpha^{p'}}\right]_1
\frac{\de}{\de z^{q'}_\beta}\right)=
\left[\frac{\de z_\beta^{q'}}{\de z_\alpha^{p'}}\right]_1
\pi\left(
\frac{\de }{\de z_\beta^{q'}}+
 \tilde\rho(\left[(a_\beta)^s_{q'}\right]_2)
\frac{\de}{\de z_\beta^s}
\right)
\\
&=&
\pi\left(
\left[\frac{\de z_\beta^{q'}}{\de z_\alpha^{p'}}\right]_2
\frac{\de}{\de z_\beta^{q'}}
+
  \left[\frac{\de z_\beta^{q'}}{\de z_\alpha^{p'}}\right]_2
\tilde{\rho}([(a_\beta)^s_{q'}]_2)
\frac{\de}{\de z_\beta^s}\right)\\
&=&
\pi\left(\frac{\de }{\de z_\alpha^{p'}}
+
\left[\frac{\de z_\beta^{q'}}{\de z_\alpha^{p'}}\right]_2
\left[\frac{\de z_\alpha^r}{\de z_\beta^{q'}}\right]_2
\frac{\de}{\de z_\alpha^r}\right)\\
& &\quad+
\pi\left( \left[\frac{\de z_\beta^{q'}}{\de z_\alpha^{p'}}\right]_2
\tilde\rho([(a_\beta)^s_{q'}]_2) \left[\frac{\de z_\alpha^r}{\de z_\beta^s}\right]_2
\frac{\de}{\de z_\alpha^r}\right)
\end{eqnarray*}
Hence
\begin{eqnarray*}
(\psi_\beta&-&\psi_\alpha)\left(\frac{\de}{\de z^{p'}_\alpha}\right)\\
&=&\pi\left(\left[\frac{\de z_\beta^{q'}}{\de z_\alpha^{p'}}\right]_1
\left[\frac{\de z_\alpha^r}{\de z_\beta^{q'}}\right]_2   \frac{\de }{\de
z_\alpha^{r}} \right)\\
& &\quad
+\pi \left(\left\{
\left[\frac{\de z_\beta^{q'}}{\de z_\alpha^{p'}}\right]_1
\tilde\rho([(a_\beta)^s_{q'}]_2)
\left[\frac{\de z_\alpha^{r}}{\de z_\beta^{s}}\right]_1
-
\tilde{\rho}([(a_\alpha)^r_{p'}]_2)
\right\}
 \frac{\de}{\de z_\alpha^{r}}\right)
\end{eqnarray*}
In Case~1 the second line of (\ref{eq:calfabeta}) yields $\psi_\beta-\psi_\alpha\equiv O$. In
Case~2, applying $\tilde\rho$ to the second line of (\ref{eq:calfabeta}), and
recalling that we are using a comfortable atlas, we get
\begin{eqnarray*}
\left[\frac{\de z_\beta^{q'}}{\de z_\alpha^{p'}}\right]_1
\tilde\rho([(a_\beta)^s_{q'}]_2)
\left[\frac{\de z_\alpha^{r}}{\de z_\beta^{s}}\right]_1
-
\tilde{\rho}([(a_\alpha)^r_{p'}]_2)
&+&
\left[\frac{\de z_\beta^{q'}}{\de z_\alpha^{p'}}\right]_1
\left[\frac{\de z_\alpha^r}{\de z_\beta^{q'}}\right]_2\\
&=&\left[\frac{\de z_\beta^{q'}}{\de z_\alpha^{p'}}\right]_1
\tilde\rho([(c_{\beta\alpha})_{q'}^{\tilde p'}]_2)[(a_\alpha)^r_{\tilde p'}]_1.
\end{eqnarray*}
Hence
\begin{equation}
(\psi_\beta-\psi_\alpha)\left(\frac{\de}{\de z^{p'}_\alpha}\right)=
[(a_\alpha)^r_{\tilde p'}]_1
\left[\frac{\de z_\beta^{q'}}{\de z_\alpha^{p'}}\right]_1
\left[\frac{\de(c_{\beta\alpha})_{q'}^{\tilde p'}}{\de z^t_\alpha}\right]_1
\pi\left([z_\alpha^t]_2\frac{\de}{\de z_\alpha^r}\right),
\label{eq:psiab}
\end{equation}
and we are done.
\end{proof}

\begin{corollary}
\label{th:psitangfol}
Let $S$ be a complex $m$-codimensional submanifold of
an $n$-dimensional complex manifold~$M$, and let $\ca F$ be a holomorphic
foliation~$\ca F$  on~$M$, of dimension $\ell\le n-m=\dim S$ tangent to~$S$. Assume that
$\ca F_S$ is a non-singular holomorphic foliation of~$S$. Then we can always find
an $\ca O_S$-morphism $\psi\colon\ca F_S\to\ca A$ such that $\theta_1\circ\psi=\id$
which is furthermore a Lie algebroid morphism.
\end{corollary}

\begin{proof}
In this case (\ref{eq:psiab}) shows that (\ref{eq:psialpha}) defines a global morphism~$\psi$.
To prove that it is a Lie algebroid morphism, it suffices to show that
$\{\tilde\pi(v_{\alpha,p'}\otimes[1]_2),\tilde\pi(v_{\alpha,q'}\otimes[1]_2)\}=O$ for all
$p'$,~$q'=m+1,\ldots,m+\ell$. But since $\ca F$ is tangent to~$S$, we have
$\tilde\pi(v_{\alpha,p'}\otimes[1]_2)=\pi(v_{\alpha,p'}\otimes[1]_2)$; hence it suffices to
show that
\begin{equation}
\{v_{\alpha,p'}\otimes[1]_2,v_{\alpha,q'}\otimes[1]_2\}\in\ca I_S\cdot\ca T^S_{M,S(1)}.
\label{eq:easydue}
\end{equation}
Now, $\{v_{\alpha,m+1},\ldots,v_{\alpha,m+\ell}\}$ is a local frame
for~$\ca F$, which is an involutive sheaf; it follows
that~$[v_{\alpha,p'},v_{\alpha,q'}]=O$, and (\ref{eq:easydue}) is an
immediate consequence.
\end{proof}

\begin{definition} Let $S$ be a complex submanifold of a complex manifold~$M$, and
$\ca F$ a holomorphic foliation of~$M$ of dimension $d\le\dim
S$. Assume that $S$ splits into~$M$, with first order
lifting~$\rho\colon\ca O_S\to\ca O_{S(1)}$ and associated
projection~$\sigma^*\colon\ca T_{M,S}\to\ca T_S$.
We shall say that {\sl $\ca F$ splits along
$\rho$} if $\gt f=O$ in $H^1\bigl(S,\Hom(\ca F^\sigma,\ca
N_S^*\otimes\ca F^\sigma)\bigr)$.
\end{definition}

\begin{corollary}
\label{th:corsplit}
Let $S$ be a complex $m$-codimensional submanifold of
an $n$-dimensional complex manifold~$M$, and let $\ca F$ be a holomorphic
foliation~$\ca F$  on~$M$, of dimension $\ell\le n-m=\dim S$ tangent to~$S$. Assume that $S$ is
comfortably embedded in~$M$ with respect to an $\ca F$-faithful first order lifting~$\rho\colon\ca
O_S\to\ca O_{S(1)}$. If
$\ca F$ splits along $\rho$ then there is an $\ca O_S$-morphism $\psi\colon\ca F^\sigma\to\ca A$
such that $\theta_1\circ\psi=\id$.
\end{corollary}

\begin{remark}
In general, the morphism $\psi$ provided by the previous corollary might not be a Lie algebroid
morphism, unless $\ell=1$.
\end{remark}

\begin{remark}
As a consequence of~(\ref{eq:yae}), we know that
\[
 \rho([(c_{\beta\alpha})^{p'}_{q'}]_1)= \rho\left(\left[\frac{\de
z_\alpha^{p'}}{\de z_\beta^{q'}}\right]_1\right) =\left[\frac{\de z_\alpha^{p'}}{\de
z_\beta^{q'}}\right]_2;
\]
therefore $\tilde\rho\bigl([(c_{\beta\alpha})^{p'}_{q'}]_2\bigr)=O$ is equivalent to
\[
[(c_{\beta\alpha})^{p'}_{q'}]_2= \left[\frac{\de z_\alpha^{p'}}{\de
z_\beta^{q'}}\right]_2;
\]
compare with Corollary~\ref{th:isoF}.
\end{remark}

Since we are using an adapted atlas, the first line of (\ref{eq:calfabeta})
yields
\begin{equation}
\tilde\rho\bigl([(c_{\beta\alpha})^{p'}_{q'}]_2\bigr)= [(a_\beta)^s_{q'}]_1
\left[\frac{\de
z_\alpha^{p'}}{\de z_\beta^s}\right]_2+\left[\frac{\de
z_\alpha^{p'}}{\de z_\beta^{q''}}\right]_1[(a_\beta)^{q''}_{q'}]_2.
\label{eq:qnove}
\end{equation}
This suggests a couple of sufficient conditions for the splitting of $\ca F$ along~$\rho$:

\begin{corollary}
\label{th:newdue}
Let $S$ be a comfortably embedded submanifold of a complex
manifold~$M$, with first order lifting~$\rho\colon\ca O_S\to\ca O_{S(1)}$ and associated
left splitting morphism~$\sigma^*\colon\ca T_{M,S}\to\ca T_S$. Let $\ca F$ be a
holomorphic foliation of~$M$ of dimension~$\ell\le\dim S$ such that $\rho$ is
$\ca F$-faithful.
Assume moreover that one of the following conditions is satisfied:
\begin{itemize}
\item[\textup{(a)}]$S$ is 2-linearizable, and $\ell=\dim S$;
\item[\textup{(b)}]$S$ is 2-linearizable, and there exists a nonsingular holomorphic foliation
of~$S$ transversal to~$\ca F^\sigma$.
\end{itemize}
Then $\ca F$ splits along $\rho$.
\end{corollary}

\begin{proof}
In both cases we can find a comfortable 2-splitting
atlas~$\{(U_\alpha,z_\alpha)\}$ adapted to~$S$, $\ca F^\sigma$
and~$\rho$ such that~$\tilde\rho([(c_{\beta\alpha})_{q'}^{p'}]_2)=0$
always. In case (a) this follows directly from (\ref{eq:qnove}),
because $m+\ell=n$; in case (b) the hypothesis implies the existence
of a comfortable 2-splitting atlas adapted to~$S$, $\ca F$
and~$\rho$ and such that $[\de z^{p'}_\alpha/\de
z_\beta^{q''}]_1=0$, and we are done.
\end{proof}

In Case 2 there is another condition ensuring the existence of the morphism~$\psi$:

\begin{lemma}
\label{th:global}
Let $S$ be a comfortably embedded submanifold of a complex
manifold~$M$, with first order lifting~$\rho\colon\ca O_S\to\ca O_{S(1)}$ and associated
left splitting morphism~$\sigma^*\colon\ca T_{M,S}\to\ca T_S$. Let $\ca F$ be a
holomorphic foliation of~$M$ of dimension~$\ell\le\dim S$ such that
$\rho$ is $\ca F$-faithful. Assume moreover that $\ca F_{S(1)}$ is (isomorphic to) the trivial
sheaf~$\ca O^{\oplus\ell}_{S(1)}$ of dimension~$\ell$ (this happens, for instance, if $\ca F$ is
globally generated by $\ell$ global vector fields). Then there exists an $\ca O_S$-morphism
$\psi\colon\ca F^\sigma\to\ca A$ such that
$\theta_1\circ\psi=\id$.
\end{lemma}

\begin{proof}
Let $v_1,\ldots,v_d$ be global generators of~$\ca F_{S(1)}$; by assumption, the $\hat
v_j=\sigma^*(v_j\otimes[1]_1)$ for $j=1,\ldots,\ell$ are global generators of~$\ca F^\sigma$.
We then define $\psi\colon\ca F^\sigma\to\ca{A}$ by setting
\[
\psi(\hat v_j)=\tilde\pi(v_j),
\]
where $\tilde\pi$ is the $\ca O_S$-morphism defined in
Proposition~\ref{th:newuno}, and then extending by $\ca
O_S$-linearity. It is then easy to check that
$\theta_1\circ\psi=\id$, and we are done.
\end{proof}

We finally have all the ingredients needed to prove our most general index theorem
for holomorphic foliations:

\begin{theorem}
\label{th:foltrans}
Let $S$ be a compact, complex, reduced, irreducible, possibly singular, subvariety of
dimension~$d$ of an $n$-dimensional complex manifold~$M$, and assume that $S$ has extendable normal
bundle. Let $\ca F$ be a (possibly singular) holomorphic
foliation~$\ca F$ on~$M$, of dimension $\ell\le d$. Assume that there exists an analytic
subset $\Sigma$ of~$S$ containing~$\bigl(\Sing(\ca F)\cap S\bigr)\cup S^{\rm sing}$ such that,
setting $S^o=S\setminus\Sigma$, we have either
\begin{itemize}
\item[\textup{(1)}] $\ca F$ is tangent to~$S^o$ and $\ca F_{S^o}$ is a non-singular holomorphic
foliation of~$S^o$ (and in this case we can take $\Sigma=\bigl(\Sing(\ca F)\cap
S\bigr)\cup S^{\rm sing}$); or
\item[\textup{(2)}]$S^o$ is
comfortably embedded in~$M$ with respect to a first order lifting~$\rho$ which is $\ca F$-faithful
outside of~$\Sigma$, and
\begin{itemize}
\item[\textup{(2.a)}]$S^o$ is 2-linearizable, and $\ell=\dim S$, or
\item[\textup{(2.b)}]$S^o$ is 2-linearizable, and there exists a nonsingular holomorphic foliation
of~$S^o$ transversal to~$\ca F^\sigma$, or
\item[\textup{(2.c)}]$\ca F_{S^o(1)}$ is (isomorphic to) the trivial
sheaf~$\ca O^{\oplus\ell}_{S^o(1)}$ of dimension~$\ell$, or, more generally,
\item[\textup{(2.d)}]$T_*(\gt f)=O$ in $H^1\bigl(S^o,\ca N_{S^o}^*\otimes \ca N_{S^o}\otimes(\ca
F^\sigma)^*\bigr)$.
\end{itemize}
\end{itemize}
Then $S$ has the Lehmann-Suwa index property of level~$\ell$ in case~$(1)$, and of
level~$\ell-\lfloor \ell/2\rfloor$ in case~$(2)$, on~$\Sigma$ with respect to~$\ca F$.
\end{theorem}

\begin{proof}
It follows from Theorem~\ref{th:generalindex},
Proposition~\ref{th:generalpsi}, Corollaries~\ref{th:psitangfol},
\ref{th:newdue}, and Lem\-ma~\ref{th:global}.
\end{proof}

\begin{remark}
Theorem~\ref{th:foltrans}.(1) is Lehmann-Suwa's theorem (see \cite{Le}, \cite{LS2} and \cite{S2});
Theorem~\ref{th:foltrans}.(2.a) generalizes both Camacho-Movasati-Sad theorem (Appendix
of \cite{CMS}) and Camacho's (\cite{C}) and Camacho-Lehmann's results (\cite{CL});
Theorems~\ref{th:foltrans}.(2.b), (2.c) and (2.d) are, as far as we know, new.
\end{remark}

\begin{example}
\label{ex:Filippo2}
We would like to compute the residue in the situation studied in~\cite{CMS}, to show that we
recover their theorem exactly. So let $S$ be a Riemann surface 2-linearizable in a complex
manifold~$M$, and let $\ca F$ be a one-dimensional foliation of~$M$ generated by a local vector
field~$v\in\ca T_{M,S}$ at a regular point~$p\in S$ which is an isolated singular point for~$\ca
F^\sigma$. If $(U_\alpha,z_\alpha)$ is a local chart at~$p$ adapted to~$S$ and to the first order
lifting~$\rho$, this means that we can write $v=a^1\de/\de z^1_\alpha+a^2\de/\de z^2_\alpha$ and
$p$ is an isolated zero of~$a_2$ on~$S$. In a pointed neighbourhood of~$p$ the
element~$v_{\alpha,2}$ defined in~(\ref{eq:newform}) is given by
\[
v_{\alpha,2}=\frac{\de}{\de z_\alpha^2}+\frac{a^1}{a^2}\frac{\de}{\de z^1_\alpha}.
\]
Therefore
\[
\psi\bigl(\sigma^*(v)\bigr)=\pi\left([a^2]_2\frac{\de}{\de z^2_\alpha}+[a^2]_2
\tilde\rho\left(\left[\frac{a^1}{a^2}\right]_2\right)\frac{\de}{\de z^1_\alpha}\right),
\]
and so
\[
\tilde X^\psi_v\left(\frac{\de}{\de z_\alpha^1}\right)=\left.a^2\frac{\de(a^1/a^2)}{\de z^1_\alpha}
\right|_S\frac{\de}{\de z^1_\alpha}.
\]
The unique (up to a constant) homogeneous symmetric polynomial of degree~1 in one variable is the
identity~$\id$; hence (\ref{eq:Grot}) yields
\[
\hbox{Res}_{\id}(\ca F;p)=\frac{1}{2\pi i}\int_{|a^2|=\varepsilon}\left.\frac{\de(a^1/a^2)}{\de
z^1_\alpha}\right|_S\, dz^2_\alpha=\hbox{Res}_p\left(\left.\frac{\de(a^1/a^2)}{\de
z^1_\alpha}\right|_S\, dz^2_\alpha\right),
\]
which is exactly the formula given in~\cite{CMS} (and used in Theorem~\ref{th:CS}).
\end{example}

\section{Holomorphic maps}
\label{sec:8}

In this final section we shall describe how to apply the strategy discussed in Section~6 using
holomorphic maps instead of foliations.

Let $S$ be an irreducible subvariety of a complex manifold~$M$, and let us denote
by~$\End(M,S)$ the space of holomorphic self-maps of~$M$ fixing~$S$ pointwise. We recall a
few definitions and facts from~\cite{ABT}.

\begin{definition}
Let $f\in{\rm End}(M,S)$, $f\not\equiv\id_M$.
The \emph{order of contact}~$\nu_f$
of~$f$ with~$S$ is defined by
\[
\nu_f=
\min_{h\in\ca O_{M,p}}\max\bigl\{\mu\in\N \mid h\circ f-h\in\ca I_{S,p}^\mu\bigr\}\in\N^*,
\]
where $p$ is any point of~$S$.
\end{definition}

In \cite{ABT}, Section~1, we proved that the order of contact is well-defined (i.e., it does not
depend on the point $p\in S$), and that it can be computed by the formula
\[
\nu_f=\min_{j=1,\ldots,n}\max\bigl\{\mu\in\N\mid f_\alpha^j-z_\alpha^j\in\ca I^\mu_{S,p}\bigr\},
\]
where $(U_\alpha,z_\alpha)$ is any local chart at~$p$, and $f_\alpha^j=z_\alpha^j\circ f$.
In particular, $[f^j_\alpha-z^j_\alpha]_{\nu_f+1}\otimes{\de/\de z^j_\alpha}$ defines a
local section of~$\ca I_S^{\nu_f}/\ca I_S^{\nu_f+1}\otimes\ca T_{M,S}=\Sym^{\nu_f}(\ca
N_S^*)\otimes\ca T_{M,S}$, that (see \cite{ABT}, Section~3) turns out to be independent of the
particular chart chosen:

\begin{definition}
  Let $f\in\End(M,S)$,
 $f\not\equiv\id_M$. The \emph{canonical section} $X_f\in
H^0\bigl(S,\Sym^{\nu_f}(\ca N_S^*)\otimes\ca T_{M,S}\bigr)$ is given by
\[
X_f=[f^j_\alpha-z^j_\alpha]_{\nu_f+1}\otimes{\de/\de z^j_\alpha}
\]
for any local chart $(U_\alpha,z_\alpha)$ at a point $p\in S$.
\end{definition}

Since we have $\Sym^{\nu_f}(\ca N_S^*)=\bigl(\Sym^{\nu_f}(\ca
N_S)\bigr)^*$, the canonical section~$X_f$ can be thought of as an $\ca
O_S$-morphism~$X_f\colon\Sym^{\nu_f}(\ca N_S)\to\ca T_{M,S}$.

\begin{definition}
Let $f\in\End(M,S)$, $f\not\equiv\id_M$. The \emph{canonical
distribution}~$\ca F_f$ associated to~$f$ (it was denoted by~$\Xi_f$ in~\cite{ABT}) is the subsheaf
of~$\ca T_{M,S}$ defined by
\[
\ca F_f=X_f\bigl(\Sym^{\nu_f}(\ca N_S)\bigr)\subseteq\ca T_{M,S}.
\]
We shall say that $f$ is \emph{tangential} if $\ca F_f\subseteq\ca T_S$.
\end{definition}

\begin{remark}
If $S$ is smooth, in \cite{ABT}, Corollary~3.2, we proved that $f$ is tangential if and only if
\[
f^r_\alpha-z^r_\alpha\in\ca I_S^{\nu_f+1}
\]
for all $r=1,\ldots,m$ and all local charts $(U_\alpha,z_\alpha)$ adapted to~$S$. We also
refer to \cite{ABT} for a discussion of the relevance of this notion.
\end{remark}

Until further notice, we shall assume that $S$ is a smooth complex submanifold of~$M$.
We shall also assume that
\begin{equation}\label{eq:limitd}
\hbox{\rm rk}_{\ca O_S}\Sym^{\nu_f}(\ca
N_S)={m+\nu_f-1\choose\nu_f}\le\dim S,
\end{equation}
where $m$ is the codimension of~$S$.

If $(U_\alpha,z_\alpha)$ is a local chart adapted to~$S$, we can find
$(g_\alpha)^j_{r_1\ldots r_{\nu_f}}\in\ca O_M$ symmetric in the lower indices such
that
\[
f^j_\alpha-z^j_\alpha=(g_\alpha)^j_{r_1\ldots r_{\nu_f}}\,z_\alpha^{r_1}\cdots
z_\alpha^{r_{\nu_f}}.
\]
The $(g_\alpha)^j_{r_1\ldots r_{\nu_f}}$ are not uniquely defined as elements of~$\ca
O_M$, but it is not difficult to check that
the $[(g_\alpha)^j_{r_1\ldots r_{\nu_f}}]_1\in\ca O_S$ \emph{are} uniquely defined.
Furthermore, the sheaf
$\ca F_f$ is locally generated by the elements
\[
v_{r_1\ldots
r_{\nu_f},\alpha}=
[(g_\alpha)^j_{r_1\ldots r_{\nu_f}}]_1\,\frac{\de}{\de z^j_\alpha}.
\]
Finally, the fact that
$X_f$ is well-defined is equivalent to the formula
\[
[(g_\alpha)^j_{r_1,\ldots,r_{\nu_f}}]_1
\left[\frac{\de z^h_\beta}{\de z^j_\alpha}\right]_1=\left[\frac{\de z_\beta^{s_1}}{\de
z_\alpha^{r_1}}
\cdots\frac{\de z_\beta^{s_{\nu_f}}}{\de z_\alpha^{r_{\nu_f}}}\right]_1
[(g_\beta)^h_{s_1,\ldots,s_{\nu_f}}]_1,
\]
so that
\begin{equation}
v_{r_1\ldots
r_{\nu_f},\alpha}=
\left[\frac{\de z_\beta^{s_1}}{\de z_\alpha^{r_1}}
\cdots\frac{\de z_\beta^{s_{\nu_f}}}{\de z_\alpha^{r_{\nu_f}}}\right]_1
v_{s_1\ldots
s_{\nu_f},\beta} .
\label{eq:mduebis}
\end{equation}

We would like to build our morphism~$\psi$ outside singularities. So in the tangential case we
shall momentarily assume that

\begin{case} The sheaf $\ca F_f$ is the sheaf of germs of holomorphic sections of a sub-bundle
of~$TS$ of rank
\[
\ell={m+\nu_f-1\choose\nu_f}.
\]
To be consistent with the non-tangential case, we shall also set
$\ca F^\sigma_f=\ca F_f$, $\sigma^*=\id_{\ca F_f}$ and $v^\sigma_{r_1\ldots r_{\nu_f},\alpha}=
v_{r_1\ldots r_{\nu_f},\alpha}$.
\end{case}

\begin{remark}
In other words, we are removing from $S$ the analytic subset of the points of $S$ where
$X$ is not injective, together with the analytic subset of the points of $S$ where $\ca T_S/\ca F_f$
is not locally free.
\end{remark}

In the non-tangential case, we project $\ca F_f$ into~$\ca T_S$, as usual.

\begin{definition} Let $S$ be a splitting submanifold of a complex manifold~$M$. Given a first order
lifting~$\rho\colon\ca O_S\to\ca O_M/\ca I_S^2$, let $\sigma^*\colon\ca T_{M,S}\to\ca T_S$
be the left splitting morphism associated to~$\rho$ by Proposition~\ref{th:ucinque}. If
$f\in\End(M,S)$, $f\not\equiv\id_M$, has order of contact~$\nu_f$, and (\ref{eq:limitd}) holds, we
shall denote by~$\ca F_f^\sigma$ the coherent sheaf of
$\ca O_S$-modules given by
\[
\ca F_f^\sigma=\sigma^*\circ X\circ(df)^{\otimes\nu_f}\bigl(\Sym^{\nu_f}(\ca
N_S)\bigr)\subseteq\ca T_S,
\]
where $(df)^{\otimes\nu_f}$ is the endomorphism of~$\Sym^{\nu_f}(\ca
N_S)$ induced by the action of~$df$ on~$N_S$. Notice that if $\nu_f>1$ (or $\nu_f=1$ and $f$ is
tangential) we have $df|_{N_S}=\id$, and hence the presence of $df$ is meaningful only
for~$\nu_f=1$ and $f$ not tangential.  We shall say that $\rho$ is \emph{$f$-faithful outside an
analytic subset~$\Sigma\subseteq S$} if
$\ca F_f^\sigma$ is the sheaf of germs of holomorphic sections of a sub-bundle of
rank~$\ell={m+\nu_f-1\choose\nu_f}$ of~$TS$ on~$S\setminus\Sigma$.  If
$\Sigma=\void$ we shall simply say that $\rho$ is \emph{$f$-faithful.}
\end{definition}

\begin{remark}
The assumption of faithfulness amounts to saying that $\sigma^*\circ X\circ(df)^{\otimes\nu_f}$ is
injective and $\ca T_S/\ca F_f^\sigma$ is locally free outside~$\Sigma$. In particular, if $m=1$
then either
$\sigma^*\circ X\circ(df)^{\otimes\nu_f}$ is identically zero or $\rho$ is $f$-faithful outside a
suitable analytic subset.
\end{remark}

So in the non-tangential case we shall assume that

\begin{case} There exists an $f$-faithful first order lifting $\rho$, with associated
left splitting morphism~$\sigma^*$. If $\{(U_\alpha,z_\alpha)\}$ is an atlas adapted to~$\rho$, we
shall also set
\[
v^\sigma_{r_1\ldots r_{\nu_f},\alpha}=\sigma^*(v_{r_1\ldots
r_{\nu_f},\alpha})=
[(g_\alpha)^p_{r_1\ldots r_{\nu_f}}]_1\,\frac{\de}{\de z^p_\alpha}
\]
when $\nu_f>1$, and
\[
v^\sigma_{r,\alpha}=\sigma^*\circ X\circ df(\de_{r,\alpha})=
\bigl[\bigl(\delta^s_r+(g_\alpha)^s_r\bigr)(g_\alpha)^p_s]_1\,\frac{\de}{\de
z^p_\alpha}
\]
when $\nu_f=1$, so that the $v^\sigma_{r_1\ldots r_{\nu_f},\alpha}$ form a local frame for~$\ca
F^\sigma_f$ and (\ref{eq:mduebis}) still holds.
\end{case}

We are now ready to compute the obstruction to the existence of the morphism~$\psi$.

\begin{proposition}
\label{th:obstmap}
Assume we are in Case~3 or in Case~4.
Given an atlas $\gt U=\{(U_\alpha,z_\alpha)\}$ adapted to~$S$
(and to $\rho$ in Case~4), let $\{m_{\beta\alpha}\}$ be the cocycle defined by
\[
m_{\beta\alpha}=
\left.\frac{\de z_\beta^{ q}}{\de z_\alpha^{p}}
\frac{\de^2 z_\alpha^{ r}}{\de z_\beta^{q} \de z_\beta^t}
\frac{\de z^t_\beta}{\de z_\alpha^s}
\,(g_\alpha)^p_{r_1\ldots r_{\nu_f}}\right|_S
\omega^s_\alpha
\otimes \de_{r,\alpha}
\otimes v^{\sigma, r_1\ldots r_{\nu_f}}_\alpha,
\]
if $\nu_f>1$ or $f$ is tangential, or by
\[
m_{\beta\alpha}=
\left.\frac{\de z_\beta^{ q}}{\de z_\alpha^{p}}
\frac{\de^2 z_\alpha^{ r}}{\de z_\beta^{q} \de z_\beta^t}
\frac{\de z^t_\beta}{\de z_\alpha^s}
\,\bigl(\delta^u_{r_1}+(g_\alpha)^u_{r_1}\bigr)(g_\alpha)^p_u\right|_S
\omega^s_\alpha
\otimes \de_{r,\alpha}
\otimes v^{\sigma, r_1}_\alpha,
\]
if $\nu_f=1$ and $f$ is not tangential, where the $v^{\sigma, r_1\ldots r_{\nu_f}}_\alpha$ form
the frame of $(\ca F^\sigma_f)^*$ dual to the frame $\{v^\sigma_{r_1\ldots r_{\nu_f},\alpha}\}$
of~$\ca F^\sigma_f$, and denote by $\gt m \in H^1\bigl(S,\ca N_S^*\otimes\ca N_S\otimes(\ca
F^\sigma_f)^*\bigr)$ the corresponding cohomology class. Then there exists a morphism
$\psi\colon{\ca F^\sigma_f}\to {\ca A}$ such that $\theta_1\circ \psi=\id$ if and
only if
\[
\gt m=O
\]
in $H^1\bigl(S,\ca N_S^*\otimes \ca N_S\otimes(\ca F^\sigma_f)^*\bigr)$.
\end{proposition}

\begin{proof}
Let us first assume $\nu_f>1$ or $f$ tangential.
We then define local $\ca O_S$-morphism $\psi_\alpha\colon\ca F^\sigma_f|_{U_\alpha}\to \ca A$
by setting
\begin{equation}
\psi_\alpha(v^\sigma_{r_1\ldots r_{\nu_f},\alpha})=
\pi\left( [(g_\alpha)^p_{r_1\ldots r_{\nu_f}}]_2
{\de \over\de z_\alpha^{p}}
\right)
\label{eq:mcinque}
\end{equation}
where $\pi\colon\ca T^S_{M,S(1)}\to\ca{A}$ is the canonical projection, and then extending
by $\ca O_S$-linearity. Notice that the argument of~$\pi$ in~(\ref{eq:mcinque}) is not well-defined,
but its image under~$\pi$ is. Since  $\theta_1\circ\psi_\alpha=\id$, it suffices to show
that the cocycle $\{m_{\beta\alpha}\}$ is represented by $\psi_\beta-\psi_\alpha$. But indeed,
recalling (\ref{eq:mduebis}) and using either that $f$ is tangential or that we are
working with an atlas adapted to~$\rho$, we have
\begin{eqnarray*}
\psi_\beta&&\!\!\!\!(v^\sigma_{r_1\ldots r_{\nu_f},\alpha}) - \psi_\beta(v^\sigma_{r_1\ldots
r_{\nu_f},\alpha})\\
&&=
\left[\frac{\de z_\beta^{s_1}}{\de z_\alpha^{r_1}}
\cdots\frac{\de z_\beta^{s_{\nu_f}}}{\de z_\alpha^{r_{\nu_f}}}\right]_1
\psi_\beta(v^\sigma_{s_1\ldots s_{\nu_f},\beta})-
 [(g_\alpha)^p_{r_1\ldots r_{\nu_f}}]_1\,\pi\left(\frac{\de}{\de z^j_\alpha}\right)
\\
&&=
[(g_\alpha)^j_{r_1\ldots r_{\nu_f}}]_1\left[\frac{\de z_\beta^q}{\de z_\alpha^j}\right]_1
\pi\left(\left[\frac{\de z_\alpha^k}{\de z_\beta^q}\right]_2\frac{\de}{\de z_\alpha^k}\right)-
 [(g_\alpha)^p_{r_1\ldots r_{\nu_f}}]_1\,\pi\left(\frac{\de}{\de z^j_\alpha}\right)
\\
&&=
[(g_\alpha)^j_{r_1\ldots r_{\nu_f}}]_1\left[\frac{\de z_\beta^q}{\de z_\alpha^j}\right]_1
\left\{\left[\frac{\de z_\alpha^p}{\de z_\beta^q}\right]_1
\pi \left(\frac{\de}{\de z^p_\alpha}\right) +
\pi\left(\left[\frac{\de z_\alpha^r}{\de z_\beta^q}\right]_2\frac{\de}{\de z^r_\alpha}
\right)\right\}\\
& &\quad- [(g_\alpha)^p_{r_1\ldots r_{\nu_f}}]_1\,\pi\left(\frac{\de}{\de z^j_\alpha}\right)
\\
&&=
[(g_\alpha)^p_{r_1\ldots r_{\nu_f}}]_1
\left[\frac{\de z_\beta^q}{\de z_\alpha^p}\right]_1
\left[\frac{\de^2 z_\alpha^r}{\de z_\beta^q\de z_\beta^t}\right]_1
\left[\frac{\de z_\beta^t}{\de z_\alpha^s}\right]_1
\pi\left([z_\alpha^s]_2\frac{\de}{\de z_\alpha^r}\right),
\end{eqnarray*}
and we are done in this case.

If $\nu_f=1$ and $f$ is not tangential we define $\psi_\alpha$ by
\[
\psi_\alpha(v^\sigma_{r,\alpha})=
\pi\left( \bigl[\bigl(\delta^s_r+(g_\alpha)^s_r\bigr)(g_\alpha)^p_s]_2
{\de \over\de z_\alpha^{p}}
\right),
\]
and the assertion follows as before.
\end{proof}

It turns out that in codimension~1 (assuming $S$ comfortably embedded in Case~4) we have $\gt m=O$
always:

\begin{proposition}
\label{th:codimuno}
Assume we are in Case 3 or in Case 4 with $S$ comfortably embedded, and that $S$ has codimension 1.
Then $\gt m=O$.
\end{proposition}

\begin{proof}
Let $\{(U_\alpha,z_\alpha)\}$ be
an atlas adapted to~$S$, and also comfortable and adapted to~$\rho$ in Case~4. We define a 0-cochain
$\{x_\alpha\}\in H^0\bigl(\gt U_S,\ca N_S^*\otimes\ca N_S\otimes(\ca F^\sigma_f)^*\bigr)$ by
setting
\begin{equation}
x_\alpha(v^\sigma_{1\ldots 1,\alpha})=\pi\left(\tilde\rho([(g_\alpha)^1_{1\ldots 1}]_2)\frac{\de}
{\de z_\alpha^1}
\right)=\pi\left(\left[\frac{\de (g_\alpha)^1_{1\ldots 1}}{\de z_\alpha^1}\right]_1[z_\alpha^1]_2
\frac{\de}{\de z_\alpha^1}\right),
\label{eq:x}
\end{equation}
where in Case~3 we have $\tilde\rho([(g_\alpha)^1_{1\ldots 1}]_2)=[(g_\alpha)^1_{1\ldots 1}]_2$,
and thus we do not need to assume anything on the embedding of~$S$ into~$M$. Notice
furthermore that, since the codimension of~$S$ is~1, the germ $(g_\alpha)^1_{1\ldots 1}$ is
well-defined \emph{as germ in~$\ca O_M$,} and not only as germ in~$\ca O_S$, and so~(\ref{eq:x}) is
well-defined.

To prove the assertion it suffices then to show that $x_\alpha-x_\beta=m_{\beta\alpha}$. Now we have
\begin{eqnarray*}
(x_\alpha&-&x_\beta)(v^\sigma_{1\ldots 1,\alpha})=
\pi\left(\tilde\rho([(g_\alpha)^1_{1\ldots 1}]_2)\frac{\de}
{\de z_\alpha^1}\right)-\left[\frac{\de z_\beta^1}{\de z_\alpha^1}\right]_1^{\nu_f}x_\beta
(v^\sigma_{1\ldots 1,\beta})\\
&=&\pi\left(\left\{\tilde\rho([(g_\alpha)^1_{1\ldots 1}]_2)-\left[\frac{\de z_\beta^1}{\de
z_\alpha^1}\right]_2^{\nu_f}\tilde\rho([(g_\beta)^1_{1\ldots 1}]_2)
\left[\frac{\de z_\alpha^1}{\de z_\beta^1}\right]_2
\right\}\frac{\de}{\de z_\alpha^1}\right).
\end{eqnarray*}
If $\nu_f>1$ from
\begin{eqnarray*}
(g_\beta)^1_{1\ldots 1}(z_\beta^1)^{\nu_f}&=&f^1_\beta-z^1_\beta=\frac{\de z_\beta^1}
{\de z_\alpha^j}(f^j_\alpha-z^j_\alpha)+R_{2\nu_f}\\
&=&
\frac{\de z_\beta^1}{\de z_\alpha^j}(g_\alpha)^j_{1\ldots 1}(z_\alpha^1)^{\nu_f}+R_{2\nu_f},
\end{eqnarray*}
where $R_{2\nu_f}\in\ca I_S^{2\nu_f}$, we get
\[
\left[\frac{\de z_\beta^1}{\de z_\alpha^1}\right]_2^{\nu_f}[(g_\beta)^1_{1\ldots 1}]_2
=\left[\frac{\de z_\beta^1}{\de z_\alpha^1}\right]_2[(g_\alpha)^1_{1\ldots 1}]_2+
\left[\frac{\de z_\beta^1}{\de z_\alpha^p}\right]_2[(g_\alpha)^p_{1\ldots 1}]_2.
\]
Since we are working with a comfortable atlas, $\tilde\rho([\de z_\beta^1/\de z_\alpha^1]_2)=O$;
therefore applying (in Case~4) $\tilde\rho$ to the previous formula we get
\[
\left[\frac{\de z_\beta^1}{\de z_\alpha^1}\right]_2^{\nu_f}\tilde\rho([(g_\beta)^1_{1\ldots 1}]_2)
=\left[\frac{\de z_\beta^1}{\de z_\alpha^1}\right]_2\tilde\rho([(g_\alpha)^1_{1\ldots 1}]_2)+
\left[\frac{\de z_\beta^1}{\de z_\alpha^p}\right]_2[(g_\alpha)^p_{1\ldots 1}]_2.
\]
Remarking that
\[
\left[\frac{\de^2 z_\beta^1}{\de z_\alpha^p\de z_\alpha^1}\right]_1
\left[\frac{\de z_\alpha^1}{\de z_\beta^1}\right]_1=
-\left[\frac{\de z_\alpha^1}{\de z_\beta^1}\right]_1
\left[\frac{\de z_\beta^q}{\de z_\alpha^p}\right]_1
\left[\frac{\de^2 z_\alpha^1}{\de z_\beta^q\de z_\beta^1}\right]_1,
\]
we obtain the assertion if $\nu_f>1$ or $f$ is tangential.

If instead $\nu_f=1$, recalling that we are using a comfortable atlas, we have
\begin{eqnarray*}
(g_\beta)^1_1 z_\beta^1&=&f^1_\beta-z^1_\beta=\frac{\de z_\beta^1}
{\de z_\alpha^j}(f^j_\alpha-z^j_\alpha)+\frac{1}{2}
\frac{\de^2 z_\beta^1}{\de z_\alpha^h\de z_\alpha^k}(f_\alpha^h-z_\alpha^h)
(f_\alpha^k-z_\alpha^k)+R_3\\
&=&
\frac{\de z_\beta^1}{\de z_\alpha^j}(g_\alpha)^j_1 z_\alpha^1+
\frac{\de^2 z_\beta^1}{\de z_\alpha^p\de
z_\alpha^1}(g_\alpha)^p_1(g_\alpha)^1_1(z_\alpha^1)^2+R_3,
\end{eqnarray*}
where $R_3\in\ca I_S^3$, and so
\[
\left[\frac{\de z_\beta^1}{\de z_\alpha^1}\right]_2[(g_\beta)^1_1]_2
=\left[\frac{\de z_\beta^1}{\de z_\alpha^1}\right]_2[(g_\alpha)^1_1]_2+
\left[\frac{\de z_\beta^1}{\de z_\alpha^p}\right]_2\bigl(1+
[(g_\alpha)^1_1]_2\bigr)[(g_\alpha)^p_1]_2.
\]
Applying $\tilde\rho$ and arguing as before we obtain the assertion
in this case too.
\end{proof}

We are now ready for our most general index theorem for holomorphic maps:

\begin{theorem}
\label{th:maphcd}
Let $S$ be a compact, complex, reduced, irreducible, possibly singular, subvariety of
codimension~$m$ of an $n$-dimensional complex manifold~$M$, and assume that $S$ has extendable
normal bundle. Let $f\in\End(M,S)$, $f\not\equiv\id_M$, have order of contact~$\nu_f$ with~$S$,
and such that $\ell={m+\nu_f-1\choose\nu_f}\le\dim S$. Assume that there
exists an analytic subset $\Sigma$ of~$S$ containing~$S^{\rm
sing}$ such that, setting $S^o=S\setminus\Sigma$, we have either
\begin{itemize}
\item[\textup{(1)}] $f$ is tangential to~$S^o$, $X|_{S^o}\colon\ca\Sym^{\nu_f}(\ca N_{S^o}^*)\to\ca
T_{S^o}$ is injective, $\ca T_{S^o}/(\ca F_f|_{S^o})$ is locally free, and
\begin{itemize}
\item[\textup{(1.a)}] $S$ has codimension~$m=1$, or, more generally,
\item[\textup{(1.b)}] $\gt m=O$;
\end{itemize}
\item[\textup{(2)}]$S^o$ is
comfortably embedded in~$M$ with respect to a first order lifting~$\rho$ which is $f$-faithful
outside of~$\Sigma$, and
\begin{itemize}
\item[\textup{(2.a)}] $S$ has codimension~$m=1$, or, more generally,
\item[\textup{(2.b)}] $\gt m=O$.
\end{itemize}
\end{itemize}
Then $S$ has the Lehmann-Suwa index property of level~$\ell-\lfloor \ell/2\rfloor$ on~$\Sigma$
with respect to~$f$.
\end{theorem}

\begin{proof}
It follows from Theorem~\ref{th:generalindex} and
Propositions~\ref{th:obstmap} and \ref{th:codimuno}.
\end{proof}

\begin{remark}
Theorems~\ref{th:maphcd}.(1.a) and (2.a) are contained
in~\cite{ABT}. Theorems~\ref{th:maphcd}.(1.b) and (2.b) are, as far
as we know, new.
\end{remark}


\begin{thebibliography}{33}

\bibitem{A}
Abate, M.: The residual index and the dynamics of holomorphic maps tangent to the
identity. Duke Math. J. {\bf 107,} 173--207 (2001).

\bibitem{ABT}
Abate, M., Bracci, F., Tovena, F.: Index theorems for
holomorphic self-maps. Ann. Math. {\bf 159,} 819--864 (2004).

\bibitem{ABT3}
Abate, M., Bracci, F., Tovena, F.: Embeddings of holomorphic varieties and normal sheaves. In
preparation (2006).

\bibitem{At}
Atiyah, M.F.: Complex analytic connections in fibre
bundles. Trans. Amer. Math. Soc. {\bf 85,} 181--207 (1957).

\bibitem{BB}
Baum, P., Bott, R.: Singularities of holomorphic foliations. J. Diff. Geometry
{\bf 7,} 279--342 (1972).

\bibitem{BS1}
Bracci, F., Suwa. T.: Residues for singular pairs and dynamics of biholomorphic maps of singular
surfaces. International J. Math. {\bf 15,} 443--466 (2004).

\bibitem{BS2}
Bracci, F., Suwa, T.: Residues for holomorphic foliations of singular pairs. Adv. Geom. {\bf 5,}
81--95 (2005).

\bibitem{BT}
Bracci, F., Tovena, F.: Residual indices of
holomorphic maps relative to  singular curves of fixed points on
surfaces. Math. Z. {\bf 242,} 481--490 (2002).

\bibitem{Bru}
Brunella, M.: Feuilletage holomorphes sue les surfaces
complexes compactes. Ann. Sci. ENS. {\bf 30,} 569--594 (1997).

\bibitem{C}
Camacho, C.: Dicritical singularities of holomorphic vector fields. Contemp.
Math. {\bf 269,} 39--45 (2001).

\bibitem{CS}
Camacho, C., Sad, P.: Invariant varieties through singularities of holomorphic
vector fields. Ann. Math. {\bf 115,} 579--595 (1982).

\bibitem{CL}
Camacho, C., Lehmann, D.: Residues of holomorphic foliations
relative to a general submanifold. Bull. London Math. Soc. {\bf 37,} 435--445 (2005).

\bibitem{CMS}
Camacho, C., Movasati, H., Sad, P.: Fibered neighborhoods of curves in surfaces.
J. Geom. Anal. {\bf 13,} 55--66 (2003).

\bibitem{CC}
Carrell, J.B., Lieberman, D.I.: Vector fields and Chern numbers. Math. Ann. {\bf 225,}
263--273 (1977).

\bibitem{Ei}
Eisenbud, D.: Commutative algebra with a view toward algebraic geometry. Springer-Verlag, New
York (1994).

\bibitem{GR}
G\'erard, R., Ramis, J.-P.: R\'esidu d'une connexion holomorphe. In: G\'erard, R., Ramis, J.-P.
(eds.),  Equations diff\'erentielles et syst\`emes de Pfaff dans le champ complexe, I\negthinspace
I, pp. 243--306. Lect. Notes in Math. 1015, Springer, Berlin (1983).

\bibitem{G}
Griffiths, P.A.: The extension problem in complex analysis I\negthinspace I; embeddings
with positive normal bundle. Amer. J. Math. {\bf 88,} 366--446 (1966).

\bibitem{Hon}
Honda, T.: Tangential index of foliations with curves on surfaces. Hokkaido Math. J.
{\bf 33,} 255--273 (2004).

\bibitem{Le}
Lehmann, D.: R\'esidus des sous-vari\'et\'es
invariants d'un feuilletage singulier. Ann. Inst. Fourier, Grenoble {\bf 41,} 211--258
(1991).

\bibitem{LS1}
Lehmann, D., Suwa, T.: Residues of holomorphic
vector fields relative to singular invariant subvarieties. J. Diff. Geom. {\bf 42,}
165--192 (1995).

\bibitem{LS2}
Lehmann, D., Suwa, T.: Generalizations of variations and Baum-Bott residues for holomorphic
foliations on singular varieties. Int. J. Math. {\bf 10,} 367--384 (1999).

\bibitem{Li}
Lins Neto, A.: Algebraic solutions of polynomial differential equations and
foliations in dimension two. In: G\'omez-Mont, X., Seade, J., Verjovski, A. (eds.),
Holomorphic dynamics, Mexico, 1986, pp. 192--232. Lect. Notes in Math. 1345, Springer, Berlin
(1988).

\bibitem{M}
Mackenzie, K.C.H.: Lie algebroids and Lie pseudoalgebras. Bull. London Math. Soc. {\bf 27,} 97--147
(1995).

\bibitem{MR}
Morrow, J., Rossi, H.: Submanifolds of $\Pro^n$ with
splitting normal bundle sequence are linear. Math. Ann. {\bf 234,} 253--261 (1978).

\bibitem{MP}
Mustata, M., Popa, M.: A new proof of a theorem of Van de
Ven. Bull. Math. Soc. Sc. Math. Roumanie {\bf 39,} 243--249 (1996).

\bibitem{PS}
Pereira, J.V., Sad, P.: On the holonomy group of algebraic curves invariant by holomorphic
foliations. Preprint series A245/2003, IMPA (2003).

\bibitem{S1}
Suwa, T.: Indices of holomorphic vector fields relative to invariant curves on surfaces.
Proc. Amer. Math. Soc. {\bf 123,} 2989--2997 (1995).

\bibitem{S2}
Suwa, T.: Indices of vector fields and residues
of singular holomorphic foliations. Hermann, Paris (1998).

\bibitem{VdV}
Van de Ven, A.: A property of algebraic varieties in
complex projective spaces. In: Colloque de g\'eom\'etrie
diff\'erentielle globale, Bruxelles, 1958, pp.~151--152. Centre Belge Rech. Math., Louvain (1959).



\end{thebibliography}
\end{document}